\documentclass[review,onefignum,onetabnum]{siamart190516}

\usepackage{lipsum}
\usepackage{amsfonts}
\usepackage{graphicx}
\usepackage{epstopdf}
\usepackage{algorithmic}
\ifpdf
  \DeclareGraphicsExtensions{.eps,.pdf,.png,.jpg}
\else
  \DeclareGraphicsExtensions{.eps}
\fi

\allowdisplaybreaks

\newcommand{\diagi}{\begin{smallmatrix}\vspace{-0.5ex}\textrm{\normalsize diag}\\\vspace{-0.8ex}i\in\mathcal{I}_{n}\end{smallmatrix}}
\newcommand{\diagik}{\begin{smallmatrix}\vspace{-0.5ex}\textrm{\normalsize diag}\\\vspace{-0.8ex}i\in\mathcal{I}_{k}\end{smallmatrix}}

\newcommand{\cred}[1]{\textcolor{red}{#1}}

\usepackage{todonotes}
\usepackage{ulem}
\usepackage{xcolor}

\newcommand{\C}{\mathbb{C}}
\newsiamremark{remark}{Remark}
\newsiamremark{hypothesis}{Hypothesis}
\crefname{hypothesis}{Hypothesis}{Hypotheses}
\newsiamthm{claim}{Claim}

\headers{Symbol-based multigrid methods for Block structures}{Matthias Bolten, Marco Donatelli, Paola Ferrari, Isabella Furci}

\title{A symbol based analysis for multigrid methods for \\ Block-Circulant and Block-Toeplitz Systems \thanks{Submitted to the editors DATE.
\funding{This work was supported by Gruppo Nazionale per il Calcolo Scientifico (GNCS-INdAM).}}}

\author{Matthias Bolten \thanks{Department of Mathematics and Informatics.
University of Wuppertal, Wuppertal, Germany 
  (\email{bolten@math.uni-wuppertal.de}, \email{furci@uni-wuppertal.de})}, \and  Marco Donatelli\thanks{Department of Science and High Technology, Insubria University,  Insubria University,  Como, Italy
  (\email{marco.donatelli@uninsubria.it},\email{pferrari@uninsubria.it})} \and Paola Ferrari\footnotemark[3] 
\and Isabella Furci\footnotemark[2]}

\nolinenumbers
\usepackage{amsopn}

\ifpdf
\hypersetup{
  pdftitle={Symbol-based multigrid methods},
  pdfauthor={Matthias Bolten, Marco Donatelli, Paola Ferrari, Isabella Furci}
}
\fi

\begin{document}

\maketitle

\begin{abstract}

In the literature, there exist several studies on symbol-based multigrid methods for the solution of linear systems having structured
coefficient matrices. In particular, the convergence analysis for such methods has been obtained in an elegant form in the case of Toeplitz matrices generated by a scalar-valued function. In the block-Toeplitz setting, that is, in the case where the matrix entries are small generic matrices instead of scalars, some algorithms have already been proposed regarding specific applications and a first rigorous convergence analysis has been performed in [M. Donatelli, P. Ferrari, I. Furci, D. Sesana, and S. Serra-Capizzano. Multigrid methods for block-circulant and block-Toeplitz large linear systems: Algorithmic proposals and two-grid optimality analysis. Numer. Linear Algebra Appl.]. However, with the existent symbol-based theoretical tools, it is still not possible to prove the convergence of many multigrid methods known in the literature. This paper aims to generalize the previous results giving more general sufficient conditions on the symbol of the grid transfer operators. In particular, we treat matrix-valued trigonometric polynomials which can be non-diagonalizable and singular at all points and we express the new conditions in terms of the eigenvectors associated with the ill-conditioned subspace. Moreover, we extend the analysis to the V-cycle method proving a linear convergence rate under stronger conditions, which resemble those given in the scalar case. In order to validate our theoretical findings, we present a classical block structured problem stemming from a FEM approximation of a second order differential problem. We focus on two multigrid strategies that use the geometric and the standard bisection grid transfer operators and we prove that both fall into the category of projectors satisfying the proposed conditions. In addition, using a tensor product argument, we provide a strategy to construct efficient V-cycle procedures in the block multilevel setting. 
\end{abstract}



\maketitle

\section{Introduction}
\label{sect:intro}

\cred{Linear systems with multilevel block-Toeplitz coefficient matrices arise in the discretization of many differential and integral equations. Among them, we mention the case of $\mathbb{Q}_{{r}} $ Lagrangian finite elements (FEM) approximation of  a second order differential problem \cite{fem} and the signal restoration problems where some of the sampling data are not available \cite{signal}. When dealing with large multilevel and multilevel block-Toeplitz systems, the performances (in terms of computational costs and iterations) of preconditioners based on circulant approximations  deteriorate \cite{negaLAA}. This is one of the many reasons why the class of multigrid methods is of great interest for the solution of such systems \cite{CCS, FS2, Chan-Sun2}. 
}

\cred{Convergence results for multigrid methods are usually based on the local Fourier analysis (LFA) \cite{LFA}, although several extensions and generalizations have been recently proposed in the literature \cite{FS1,H90,Nap11}. In \cite{D10} it was proven that the convergence analysis of multigrid methods for circulant and Toeplitz matrices \cite{NM,ST} is a linear algebra generalization of the LFA in the case of the Galerkin approach. Indeed, it does not necessarily require a differential operator and it can also be applied to integral problems with applications such as signal and image processing \cite{tik}. For differential problems with constant coefficients and uniform grids, the matrix algebra approach leads to a condition on the symbols associated to the circulant matrices analogous to the classical condition on the orders of the grid transfer operators.  In particular, such condition relates the order of the zeros of the symbols associated to the coefficient matrix and the grid transfer operators.} 

\cred{In this paper, we prove a generalization of such condition to block-symbols, that is when the generating function $\mathbf{f}$ associated with the coefficient matrix is a matrix-valued trigonometric polynomial. The block-symbol has been previously investigated in the literature \cite{DFFSS,FRTBS, Huck}, but many theoretical aspects have not yet been properly addressed. In particular, a V-cycle convergence analysis is still missing and some classical grid transfer operators do not satisfy the strong requirements of the two-grid analysis in \cite{DFFSS}. 
The main aim of the paper is to provide a complete convergence analysis of multigrid methods for structured block-Toeplitz and circulant systems under weak assumptions. In order to show the applicability of our theory, we will consider classical multigrid strategies for $\mathbb{Q}_{{r}} $ Lagrangian FEM in the case of uniform Cartesian grids.}

First, we consider the Two Grid Method (TGM) and, according to the classical  Ruge-St\"uben  \cite{RStub} convergence analysis, we focus on validating both a smoothing property and an approximation property. The first is easily generalizable in the block setting from the scalar case. Indeed, in \cite{DFFSS} it has been proven that it mainly affects the choice of the specific relaxation parameter for the selected smoother. The validation of the approximation property for block structured matrices approximation, instead,  is non-trivial and it requires additional hypotheses. In particular, the idea is to focus on the crucial choice of conditions on the trigonometric polynomial $\mathbf{
p}$ used to construct the projector. In   \cite{DFFSS} the proof of the approximation property is based on the validation of an additional commutativity requirement \cite[Section 4.1]{DFFSS}. However, in some practical cases these conditions cannot be satisfied. Hence, the main \textcolor{red}{theorem} of Section \ref{sec:Theorem} provides less restrictive conditions on $\mathbf{
p}$. Indeed, differently from \cite{DFFSS}, the conditions are expressed in terms of the eigenvectors associated with the
ill-conditioned subspace and permit to enlarge the class of suitable trigonometric polynomials used to
construct the projectors. Moreover, we provide some useful lemmas which  can further simplify the validation of the requirements, under specific hypotheses  that $\mathbf{p}$ often satisfies in the applications.

Another important result of the paper concerns the extension of the theoretical findings to V-cycle methods. Indeed,  following the proof of the main theorem on the TGM convergence and the results in \cite{NN}, it is possible to obtain elegant conditions for the convergence and optimality of the V-cycle in the block case. For the latter, a crucial point is the investigation of the properties of the symbols at coarse levels, with a particular focus on the orders of the zeros.

In order to validate our theoretical findings and show their applicability, we present a classical block structured problem stemming from the $\mathbb{Q}_{\textcolor{black}{r}} $ Lagrangian FEM approximation of  a second order differential problem. We focus on two multigrid strategies that use the geometric projection operator and the standard bisection grid transfer operator. We prove that both fall into the category of projectors  satisfying the proposed conditions, which lead to convergence and optimality of multigrid methods  \cite{Hack, FRTBS}. Finally, in  Section  \ref{sect:multiD} we provide the extension of the theory for block multilevel Toeplitz matrices exploiting the properties of the Kronecker product.

The paper is organized as follows. In Section \ref{sect:two_grid} we recall the basics of the multigrid methods, with particular attention to the TGM convergence analysis and on the general conditions that lead to the V-cycle optimality. 
 In Section  \ref{sec:proposal_block_circ} we restrict the attention to the block setting recalling some properties of block-circulant and block-Toeplitz matrices. In particular,  we introduce the main ingredients for an effective multigrid procedure that will be investigated in Section \ref{sec:Theorem}. Here, we focus on the conditions which ensure the convergence and optimality of the TGM for a linear system with coefficient matrix generated by a matrix-valued trigonometric polynomial and we provide a possible simplification for the validation of the conditions in practical cases. In Subsection \ref{sec:condition_Vcycle} we derive the conditions for convergence and optimality also for the V-cycle. In Section \ref{sect:natural_projectors}  we present the two classical multigrid strategies using the geometric projection operator and the standard bisection grid transfer operator. Finally, in Section \ref{sect:multiD} we consider the extension in the block multilevel case and we show how the results of previous sections can be exploited and generalized.

\section{Multigrid methods}\label{sect:two_grid}
Multigrid methods are efficient and robust iterative methods for solving linear systems of the form
\[A_{n}x_{n}=b_{n},\]
where often, and as assumed in this paper, $A_{n}\in\mathbb{C}^{n\times n}$ is positive definite \cred{\cite{Trot}}.
The main idea is to combine a classical stationary iterative method, called smoother, with a coarse grid correction having a spectral behaviour complementary with respect to the smoother \cred{\cite{tutorial}}.
In order to fix the notation for positive definite matrices, if $X\in\mathbb{C}^{n\times n}$ is a positive definite matrix,
$\|\cdot\|_{X}=\|X^{1/2}\cdot\|_{2}$ denotes the Euclidean norm
weighted by $X$ on $\mathbb{C}^{n}$. 
If $X$ and $Y$ are Hermitian matrices, then
the notation $X\leq Y$ means that $Y-X$ is a nonnegative definite matrix. \textcolor{red}{Given a matrix $X$ we denote by $X^T$ and $X^H$ the transpose and the conjugate transpose matrix of $X$, respectively.  }

\subsection{Two-grid method}
Let $P_{n,k}\in\mathbb{C}^{n\times k}$, $k<n$, be a given full-rank
matrix and let us consider two stationary iterative methods: the method $\mathcal{V}_{n,\rm{pre}}$, with iteration matrix ${V}_{n,\rm{pre}}$, and $\mathcal{V}_{n,\rm{post}}$, with iteration matrix ${V}_{n,\rm{post}}$.
An iteration of a Two-Grid Method (TGM) is given in Algorithm 1.

\begin{small}
\begin{algorithm}
	\caption{TGM$(A_{n},\mathcal{V}_{n,\rm{pre}}^{\nu_{\rm{pre}}},\mathcal{V}_{n,\rm{post}}^{\nu_{\rm{post}}},P_{n,k},b_{n},x_{n}^{(j)})$}
	\label{alg:Two_grid}
	\begin{algorithmic}
		\STATE{ 0. $\tilde{x}_{n}=\mathcal{V}_{n,\rm{pre}}^{\nu_{\rm{pre}}}(A_{n},{b}_{n},x_{n}^{(j)})$}
		\STATE{ 1. $r_{n}=b_{n}-A_{n}\tilde{x}_{n}$}
		\STATE{ 2. $r_{k}=P_{n,k}^{H}r_{n}$}
		\STATE{ 3. $A_{k}=P_{n,k}^{H}A_{n}P_{n,k}$}
		\STATE{ 4. Solve $A_{k}y_{k}=r_{k}$}
		\STATE{ 5. $\hat{x}_{n}=\tilde{x}_{n}+P_{n,k}y_{k}$}
		\STATE{ 6. $x_{n}^{(j+1)}=\mathcal{V}_{n,\rm{post}}^{\nu_{\rm{post}}}(A_{n},{b}_{n},\hat{x}_{n})$}
		
	\end{algorithmic}
\end{algorithm}
\end{small}

Steps $1.\rightarrow5.$ define the ``coarse grid correction'' that depends on the projecting operator $P_{n,k}$, while step $0.$ and
step $6.$ consist, respectively, in applying $\nu_{\rm{pre}}$ times a pre-smoother
and $\nu_{\rm{post}}$ times a post-smoother of the given iterative methods.
Step 3. defines the coarser matrix $A_{k}$ according to the Galerkin approach which ensures that the coarse grid correction is an algebraic projector and hence is very useful for an algebraic study of the convergence of the method. Indeed, the TGM is a stationary method defined by the following iteration matrix
\begin{small}
\begin{align*}
{\rm TGM}(A_{n},V_{n,\rm{pre}}^{\nu_{\rm{pre}}},V_{n,\rm{post}}^{\nu_{\rm{post}}},P_{n,k})=
V_{n,\rm{post}}^{\nu_{\rm{post}}}
\left[I_{n}-P_{n,k}\left(P_{n,k}^{H}
A_{n}P_{n,k}\right)^{-1}P_{n,k}^{H}A_{n}\right]V_{n,\rm{pre}}^{\nu_{\rm{pre}}}.
\end{align*}
\end{small}
\cred{
\begin{theorem}\label{teoconv}(\cite{RStub})
	Let $A_{n}$ be a positive definite matrix of size $n$ and let $V_{n,{\rm post}},$ $V_{n,{\rm pre}}$ be defined as in the {\rm TGM} algorithm.
	Assume
	\begin{itemize}
	\begin{small}
		\item[(a)] $\exists\alpha_{\rm{pre}}>0\,:\;\|V_{n,\rm{pre}}x_{n}\|_{A_{n}}^{2}\leq\|x_{n}\|_{A_{n}}^{2}-\alpha_{\rm{pre}}\|V_{n,\rm{pre}}x_{n}\|_{A_{n}^2}^{2},\qquad \forall x_{n}\in\mathbb{C}^{n},$
		\item[(b)] $\exists\alpha_{\rm{post}}>0\,:\;\|V_{n,\rm{post}}x_{n}\|_{A_{n}}^{2}\leq\|x_{n}\|_{A_{n}}^{2}-\alpha_{\rm{post}}\|x_{n}\|_{A_{n}^2}^{2},\qquad \forall x_{n}\in\mathbb{C}^{n},$
		\item[(c)] $\exists\gamma>0\,:\;\min_{y\in\mathbb{C}^{k}}\|x_{n}-P_{n,k}y\|_{2}^{2}\leq \gamma\|x_{n}\|_{A_{n}}^{2},\qquad \forall x_{n}\in\mathbb{C}^{n}.$
	\end{small}
	\end{itemize}
	Then $\gamma\geq\alpha_{\rm{post}}$ and
\begin{small}
	\begin{align*}
	\|{\rm TGM}(A_{n},V_{n,\rm{pre}},V_{n,\rm{post}},P_{n,k})\|_{A_{n}}\leq\sqrt{\frac{1-\alpha_{\rm{post}}/ \gamma}{1+\alpha_{\rm{pre}}/ \gamma}}<1.
	\end{align*}
\end{small}
\end{theorem}}
Conditions \cred{$(a)-(b)$} and $(c)$ are usually called ``smoothing property'' and
``approximation property'', respectively.
Since $\alpha_{\rm{post}}$ and $\gamma$ are independent of $n$, if the assumptions of Theorem \ref{teoconv} are satisfied, then the resulting TGM exhibits a linear convergence. In other words, the number of iterations in order to reach a given accuracy $\epsilon$ can be bounded from above by a constant independent of $n$ (possibly depending on the parameter $\epsilon$). Moreover, if the projection and smoothing steps have a computational cost lower or equal to the matrix-vector product with the matrix $A_n$, then the TGM is optimal.


\subsection{V-cycle method}
For large $n$ a V-cycle method should be implemented. The standard V-cycle method is obtained replacing the direct solution at step 4. with a recursive call of the TGM applied to the coarser linear system $A_{k_{\ell}}y_{k_{\ell}}=r_{k_{\ell}}$, where $\ell$ represents the level. The recursion is usually stopped at level ${\ell_{\min}}$ when $k_{\ell_{\min}}$ becomes small enough for solving cheaply step 4. with a direct solver. 
In the following, \textcolor{red}{as it is done in \cite{NN}, we assume that we are using the same iterative method as pre/post smoother, with same numbers of iterative steps. We denote the iteration matrix by $V_{n_{\ell}}^{\nu}$.}
The global iteration matrix ${\rm MGM}_{0}$ of the V-cycle method is recursively defined as
\begin{small}
\begin{align*}
&{\rm MGM}_{\ell_{\min}}(A_{n_{\ell_{\min}}},V_{n_{\ell_{\min}}}^{\nu},V_{n_{\ell_{\min}}}^{\nu},P_{n_{\ell_{\min}},k_{\ell_{\min}}})
=O_{n_{\ell_{\min}},n_{\ell_{\min}}},\\
&{\rm MGM}_\ell(A_{n_{\ell}},V_{n_{\ell}}^{\nu},V_{n_{\ell}}^{\nu},P_{n_{\ell},k_{\ell}})=\\
&
 V_{n_{\ell}}^{\nu}
\left[I_{n_{\ell}}-P_{n_{\ell},k_{\ell}}\left(I_{n_{\ell}+1}-{\rm MGM}_{\ell+1}\right)\left(P_{n_{\ell},k_{\ell}}^{H}
A_{n_{\ell}}P_{n_{\ell},k_{\ell}}\right)^{-1}P_{n_{\ell},k_{\ell}}^{H}A_{n_{\ell}}\right]V_{n_{\ell}}^{\nu}, 
\end{align*}
\end{small}
for $ \ell={\ell_{\min}-1}, \dots, 0$, \textcolor{red}{where $O_{n_{\ell_{\min}},n_{\ell_{\min}}}$ denotes the $n_{\ell_{\min}}\times n_{\ell_{\min}}$ matrix of all zero components.}

In order to prove the convergence and optimality of the V-cycle method, the key ingredient is the analysis of the spectral radius of ${\rm MGM}_0$, which is the iteration matrix at the finest level.
In \cite[Corollary 3.1]{NN} the authors show the following relation
\begin{small}
\begin{equation*}
\begin{split}
 \rho\left({\rm MGM_0}(A_{n_0},V_{n_0}^{\nu},V_{n_0}^{\nu},P_{n_0,k_0})\right)\le  1-\min_\ell \frac{1-\rho\left({\rm TGM}(A_{n_0},V_{n_0}^{\nu},V_{n_0}^{\nu},P_{n_0,k_0})\right)}{\|\pi_{A_{n_\ell}}\|^2_{A_{n_\ell}(I_{n_\ell}-V_{n_\ell}^{2\nu})^{-1}},
}
\end{split},
\end{equation*}
\end{small}
where $\pi_{A_{n_\ell}} = P_{n_\ell,k_\ell}(P_{n_\ell,k_\ell}^H A_{n_\ell}P_{n_\ell,k_\ell})^{-1}P_{n_\ell,k_\ell}A_{n_\ell}.$
Hence, assuming that we choose smoothers and prolongation operators such that the two grid optimality is guaranteed, i.e. $\rho\left({\rm TGM}(A_{n_0},V_{n_0}^{\nu},V_{n_0}^{\nu},P_{n_0,k_0})\right)<c<1$, it is sufficient to prove that the following quantity is bounded
\begin{small}
\begin{equation}\label{eq:bound_notay_pi}
\|\pi_{A_{n_\ell}}\|^2_{A_{n_\ell}(I_{n_\ell}-V_{n_\ell}^{2\nu})^{-1}}.
\end{equation}
\end{small}
In practice, the boundness of $\|\pi_{A_{n_\ell}}\|^2_2$ is usually enough, reducing the convergence analysis to the study of the spectral behaviour of the coarse grid correction operator.
\begin{lemma}\label{lem:V}
Assume that there exists a positive $C$ independent of $n$ such that
\begin{small}
\begin{equation}\label{eqn:vcycle_norm_f}
\Lambda(A_{n_\ell}) \subseteq (0,C],
\end{equation}
\end{small}
\textcolor{red}{where $\Lambda(A_{n_\ell})$ denotes the spectrum of the matrix $A_{n_\ell}$. Suppose that } one iteration of the Richardson method with the damping parameter $\omega \in(0, 2/C)$ both as pre-smoother and post-smoother is applied. Then, the boundness of $\|\pi_{A_{n_\ell}}\|^2_2$  implies that $\|\pi_{A_{n_\ell}}\|^2_{A_{n_\ell}(I_{n_\ell}-V_{n_\ell}^{2\nu})^{-1}}$ is bounded as well.
\begin{proof}
Applying one step of pre-smoother and post-smoother, it holds 
\begin{small}
\begin{eqnarray}
\|\pi_{A_{n_\ell}}\|^2_{A_{n_\ell}(I_{n_\ell}-V_{n_\ell}^{2})^{-1}}
&=&\left\|A_{n_\ell}^{1/2}(I_{n_\ell}-V_{n_\ell}^{2})^{-1/2}\pi_{A_{n_\ell}} \left(A_{n_\ell}^{1/2}(I_{n_\ell}-V_{n_\ell}^{2})^{-1/2}\right)^{-1}   \right\|_2^{2} \nonumber \\
& \le&\left\|A_{n_\ell}^{1/2}(I_{n_\ell}-V_{n_\ell}^{2})^{-1/2}\right\|_2^2\left\|\pi_{A_{n_\ell}}\right\|_2^2 \left\|\left(A_{n_\ell}^{1/2}(I_{n_\ell}-V_{n_\ell}^{2})^{-1/2}\right)^{-1}   \right\|_2^{2}. \label{eq:boundPi}
\end{eqnarray} 
\end{small}
Since the smoother is the Richardson method with damping parameter $\omega$, we can write
\begin{small}
\begin{equation*}
\begin{split}
&A_{n_\ell}^{1/2}(I_{n_\ell}-V_{n_\ell}^{2})^{-1/2}=
A_{n_\ell}^{1/2}(I_{n_\ell}-\left(I_{n_\ell}-\omega A_{n_\ell}\right)^2)^{-1/2}=
(2\omega I_{n_\ell}-\omega^2A_{n_\ell})^{-1/2},
\end{split}
\end{equation*}
\end{small}
whose 2-norm is bounded if $\omega$ is such that $
2\omega -\omega^2 \lambda_j(A_{n_\ell})>0, \, j=1,\dots,n_\ell.$
Equivalently, 
\begin{small}
\begin{equation*}
\begin{split}
&\left(A_{n_\ell}^{1/2}(I_{n_\ell}-V_{n_\ell}^{2})^{-1/2}\right)^{-1}=
(2\omega I_{n_\ell}-\omega^2A_{n_\ell})^{1/2},
\end{split}
\end{equation*}
\end{small}
whose 2-norm is bounded for all $\omega \in (0, 2/C)$ as $n$ increases thanks to equation (\ref{eqn:vcycle_norm_f}).
Finally, thanks to inequality \eqref{eq:boundPi}, the boundness of $\|\pi_{A_{n_\ell}}\|^2_2$  implies that {\small $$\|\pi_{A_{n_\ell}}\|^2_{A_{n_\ell}(I_{n_\ell}-V_{n_\ell}^{2\nu})^{-1}}$$ }is bounded as well.
\end{proof}
\end{lemma}

 \section{Multigrid methods for block-circulant and block-Toeplitz matrices}\label{sec:proposal_block_circ}
In the present paper, we are interested in proposing an effective multigrid method in the case where $A_{n}$ is a block-circulant or block-Toeplitz matrix. Therefore, we recall some properties of these structured matrices.

\subsection{Block-circulant and block-Toeplitz matrices} 
Let $\mathcal{M}_d$
be the linear space of the complex $d\times d$ matrices. Given a function $\mathbf{f}:Q\to\mathcal M_d$,  for $i=1,\dots, d$, we denote by  $\lambda_i(\mathbf{f})$ the eigenvalue functions of $\mathbf{f}$ and by $\lambda_i(\mathbf{f}(\theta))$ their evaluation at a point $\theta\in Q$.
The following lemma is derived from the results in \cite[Section VI.1]{Bhatia} and provides the existence and continuity of the eigenvalue functions of $\mathbf{f}$.
\begin{lemma}\label{lemm:continuity_eig}
Let $\theta\rightarrow \mathbf{f}(\theta)$ be a continuous map from an interval $Q$ into the space of $d\times d$ matrices such that the eigenvalues of $\mathbf{f}(\theta)$ are real for all $\theta \in Q$. Then there exist continuous functions $\lambda_1(\mathbf{f}(\theta)), \lambda_2(\mathbf{f}(\theta)), \dots, \lambda_d(\mathbf{f}(\theta))$ that, for each $\theta \in Q$, are the eigenvalues of $\mathbf{f}(\theta)$.  
\end{lemma} 

Let $\mathbf{f}:Q\to\mathcal M_d$, with $Q=(-\pi,\pi)$. We say that $\mathbf{f} \in L^p([-\pi,\pi])$ (resp. is
measurable) if all its components $f_{ij}:Q\to\mathbb C,\
i,j=1,\ldots,d,$ belong to $L^p([-\pi,\pi])$ (resp. are measurable) for $1\le
p\le\infty$. 
\begin{definition}\label{def-multilevel}
	Let the Fourier coefficients of  a function $\mathbf{f} \in L^p([-\pi,\pi])$ be
	\begin{align*}
	\hat{\mathbf{f}_j}:=\frac1{2\pi}\int_{Q}\mathbf{f}(\theta){\rm e}^{- \iota j \theta} d\theta\in\mathcal M_d,
	\qquad  \iota^2=-1, \, j\in\mathbb Z.
	\end{align*}
	Then, the block-Toeplitz matrix associated with \textcolor{red}{$\textbf{f}$ is the matrix with $d$ blocks of size $n$ and hence it has order $d\cdot n$} given by
\begin{small}
	\begin{align*}
	T_n(\mathbf{f})=\sum_{|j|<n}J_{n}^{(j)}\otimes\hat{\mathbf{f}_j},
	\end{align*}
\end{small}
	where $\otimes$ denotes the (Kronecker) tensor product of matrices. The term
	$J_n^{(j)}$ is the matrix of order $n$ whose $(i,k)$ entry equals $1$ if $i-k=j$
	and zero otherwise.  
\end{definition}
The set $\{T_n(\mathbf{f})\}_{n\in\mathbb N}$ is
called the \textit{family of block-Toeplitz matrices generated by $\mathbf{f}$}, that
in turn is referred to as the \textit{generating function or the symbol of
	$\{T_n(\mathbf{f})\}_{n\in\mathbb N}$}.
In the scalar case, when $d=1$, if ${f}$ is a trigonometric polynomial of degree lower than $n$, then we can define the circulant matrix  generated by $f$ by
\begin{small}
\begin{align*}
\mathcal{A}_{n}(f)=F_n \diagi({f}(\theta_i^{(n)}))F_{n}^{H}, 
\end{align*} 
\end{small}
where $\mathcal{I}_{n} = \{0,\ldots,n-1\}$ and $F_{n}=\frac{1}{\sqrt{n}}\left[{\rm e}^{-\imath j\theta_i^{(n)}}\right]_{i,j=0}^{n-1}$, with $\theta_i^{(n)}=\frac{2\pi i}{n}$. Circulant matrices form an algebra of normal matrices.
In the block-case, $d>1$, if $\mathbf{f}:Q\to\mathcal M_d$ is a $d\times d$ matrix-valued trigonometric polynomial, then the block-circulant matrix of order $dn$ generated by $\mathbf{f}$ is defined as
\begin{small}
  \begin{equation*}
\mathcal{A}_{n}(\mathbf{f})=(F_n\otimes I_{d}) \diagi(\mathbf{f}(\theta_i^{(n)}))(F_{n}^{H}\otimes I_{d}),
\end{equation*} 
  \end{small}  
where $\diagi(\mathbf{f}(\theta_i^{(n)}))$ is the block-diagonal matrix where the block-diagonal elements are  $\mathbf{f}(\theta_i^{(n)})$.

\subsection{Projectors for block structured matrices}\label{sect:project_circ}
For the convergence analysis of block-circulant and block-Toeplitz matrices, previous results are based on the Ruge-St\"uben theory  \cite{RStub} for TGM in Theorem \ref{teoconv}, see \cite{CCS, HS2, DFFSS}. 
The smoothing property is satisfied by damped Richardson iteration simply choosing the damping parameter in the interval $(0,2/\|\mathbf{f}\|_\infty)$, see \cite[Lemma 1]{DFFSS}.
The approximation property $(c)$ requires a precise definition of $P_{n,k}$ and a detailed analysis.

The choice of  prolongation and restriction operators  fulfilling the approximation condition is crucial for multigrid convergence and optimality.
In particular, the projector  $P_{n,k}$ is chosen is order that
\begin{itemize}
\item it projects the problem onto a coarser space by ``cutting'' the coefficient matrix,
\item the resulting projected matrix should maintain the same block structure and properties of the original matrix.
\end{itemize}

Let $K_{n,k}$ be the $n\times k$ downsampling matrix, such that:
\begin{description}
\item[n even:] $k=\frac{n}{2}$ and $K_{n,k}=K^{Odd}_{n,k}$,
\item[n odd:] $k=\frac{n-1}{2}$ and  $K_{n,k}=K^{Even}_{n,k}$,
\end{description}
with $K^{Odd}_{n,k}$ and $K^{Even}_{n,k}$ defined as
\begin{small}
\[
K^{Odd}_{n,k} = \left[\begin{array}{cccccccc}
		1 &  & & &\\
		0 &  & & &\\
   \vdots & 1& & &\\
		  & 0& & &\\	
		  &\vdots& & &\vdots\\
		  &  & & &1\\
	      &  & & &0\\		
		\end{array}\right]_{n\times k},
\qquad
K^{Even}_{n,k} = \left[\begin{array}{cccccccc}
		0 &  & & &\\
		1 &  & & &\\
        0 & 0& & &\\
	\vdots& 1& & &\\	
		  & 0& & &\\
		  & \vdots& & &\\
		  &  & & &\vdots\\
		  &  & & &1\\
	      &  & & &0\\		
		\end{array}\right]_{n\times k}.
\]
\end{small}
In particular, $K^{Odd}_{n,k}$ is the $n\times k$ matrix obtained by removing the even rows from the identity matrix of size $n$, that is it keeps the odd rows. On the other hand, $K^{Even}_{n,k}$ keeps the even rows.
When $n$ is even, $K^{Odd}_{n,k}$ performs the packaging of the Fourier frequencies since it holds
{\small \[(K^{Odd}_{n,k})^TF_{n}=\frac{1}{\sqrt{2}}[F_{k} \, | \, F_{k}].\]}
This property of the Fourier matrix is the key to define a projector $P_{n,k}$ that preserves the block-circulant structure at the coarser levels.
In the rest of the paper, the projector will be denoted by $P^d_{n,k}$ since the block-structured matrices have blocks of order $d$.
Therefore, we define the structure of the projecting operators $P^d_{n,k}$ for the block-circulant matrix $\mathcal{A}_{n}(\mathbf{f})$ generated by a trigonometric polynomial $\mathbf{f}:Q\rightarrow\mathcal{M}_d$ as follows.
Let $n$ be even and of the form $2^t$, $t\in \mathbb{N}$,  such that the size of the coarser problem is $k=\frac{n}{2}=2^{t-1}$.
The projector $P^d_{n,k}$ is then constructed as the product between a matrix $\mathcal{A}_{n}(\mathbf{p})$ in the algebra, with $\mathbf{p}$ a proper trigonometric polynomial that will be defined in the following sections, and a cutting matrix $K^{Odd}_{n,k}\otimes I_d$. That is,
\begin{small}
\begin{equation}\label{eqn:def_Pnkd}
	P^d_{n,k}=\mathcal{A}_{n}(\mathbf{p})(K^{Odd}_{n,k} \otimes I_d).
\end{equation}
\end{small}
The result of multiplying a $d\times d$ block matrix of dimension $dn\times dn$ by $K^{Odd}_{n,k}\otimes I_d$  is a $d\times d$ block matrix where just the even ``block-columns'' are maintained.
{ We are left to determine the conditions to be satisfied
by $\mathcal{A}_{n}(\mathbf{p})$ (or better by its generating function $\mathbf{p}$), in order to  obtain a projector which is effective in terms of convergence.
\cred{Using the block-symbol analysis}, sufficient conditions have been proven in \cite{DFFSS}. Unfortunately, such conditions are quite strong and are not satisfied by the classical projector studied in Section \ref{sect:natural_projectors}. Therefore, in the next section we prove weakly conditions on $\mathbf{p}$ that provide an optimal multigrid method.
}
The same strategy can be applied when we deal with block-Toeplitz matrices generated by a matrix-valued trigonometric polynomial, instead of block-circulant matrices. 
Indeed, the only thing that should be adapted is the  structure of the projector which slightly changes for block-Toeplitz matrices, in order to preserve the structure at coarser levels. 
Hence, for a matrix-valued trigonometric polynomial $\mathbf{p}$, the projector matrix is $
P^d_{n,k}=T_{n}(\mathbf{p})\left(K^{Even}_{n,k}\otimes I_d\right).$
Note that in the Toeplitz case $n$ should be chosen odd and of the form $2^{t}-1$, $t\in \mathbb{N}$,  such that the size of the coarser problem is $k=\frac{n-1}{2}=2^{t-1}-1$.

\section{Multigrid convergence for block-Circulant matrices} \label{sec:Theorem}
Let $A_{N}=\mathcal{A}_{n}(\mathbf{f})$, $N = N(d,n) = dn$ and $n$ even, with $\mathbf{f}$ matrix-valued trigonometric polynomial $\mathbf{f}\geq0$.
\textcolor{red}{We highlight that the theoretical results we derive are based on the hypothesis  that \textbf{f} is a trigonometric polynomial. This guarantees that the symbol at the coarse levels maintains the same structure and properties of the symbol at the finest one. However, in line with the scalar-valued case addressed in \cite{NM,ST}, the proposed theory can be easily extended to the dense case where \textbf{f} belongs to $ L^\infty([-\pi,\pi])$, only requiring the additional hypothesis that \textbf{f} has isolated zeros of finite order.}

Let $P^d_{n,k}=\mathcal{A}_{n}(\mathbf{p})(K^{Odd}_{n,k} \otimes I_d)$ with $\mathbf{p}$ matrix-valued trigonometric polynomial. Suppose that there exist unique $\theta_0\in[0,2\pi)$ and $\bar{\jmath} \in\{1,\dots,d\}$ such that 
\begin{small}
\begin{equation}\label{eqn:condition_on_f}
\left\{\begin{array}{ll}
\lambda_j(\mathbf{f}(\theta))=0, & \mbox{for } \theta=\theta_0 \mbox{ and } j=\bar{\jmath}, \\
\lambda_j(\mathbf{f}(\theta))>0, & {\rm otherwise}.
\end{array}\right.
\end{equation}
\end{small}
The latter assumption means that the matrix $\mathbf{f}(\theta)$ has exactly one zero eigenvalue in $\theta_0$ and it is positive definite in $[0,2\pi)\backslash\{\theta_0\}$. Moreover, we have that the order of the zero in $\theta_0$ must be even. As a consequence, the matrices $A_{N}$ could be singular and the ill-conditioned subspace is the eigenspace associated with $\lambda_{\bar{\jmath}}(\mathbf{f}(\theta_0))$. On the other hand, the block-Toeplitz matrices $T_{n}(\mathbf{f})$ are positive definite with the same ill-conditioned subspace and become ill-conditioned as $N$ increases. Since $\mathbf{f}(\theta)$ is Hermitian, it can be diagonalized by an orthogonal matrix $Q(\theta)$.  { Moreover, we are in the setting that the eigenvalues and the eigenvectors of $\mathbf{f}$ are continuous functions in the variable $\theta$ \cite{Kato, Re}. } We have
\begin{small}
\begin{equation}\label{eq:dec_f}
\begin{split}
& \mathbf{f}(\theta) = Q(\theta)D(\theta)Q(\theta)^H =\\
& \left[
\begin{array}{@{\;}c@{\;}|@{\;}c@{\;}|@{\;}c@{\;}|@{\;}c@{\;}|@{\;}c@{\;}|@{\;}c@{\;}}
q_1(\theta)&\dots &q_{\bar{\jmath}}(\theta)&\dots& q_d(\theta)
\end{array}
\right]\begin{bmatrix}
\lambda_1(\mathbf{f}(\theta))& & &  & &\\
& \ddots & &  &\\
& & \lambda_{\bar{\jmath}}(\mathbf{f}(\theta))& & \\
& & & \ddots & \\
& & & &\lambda_d(\mathbf{f}(\theta))
\end{bmatrix}
\left[
\begin{array}{ccccccc}
{q_1}^H(\theta) \\
\hline
\vdots\\
\hline
q_{\bar{\jmath}}^H(\theta)\\
\hline
\vdots\\
\hline
{q_d}^H(\theta)
\end{array}
\right],
\end{split}
\end{equation}
\end{small}
where $q_{\bar{\jmath}}(\theta)$ is the eigenvector that generates the ill-conditioned subspace since $q_{\bar{\jmath}}(\theta_0)$ is the eigenvector of $\mathbf{f}(\theta_0)$ associated with $\lambda_{\bar{\jmath}}(\mathbf{f}(\theta_0))=0$. 
Under the following  assumptions, we will prove that there are sufficient conditions to ensure the linear convergence of the TGM. 

\textcolor{red}{In the next section---in particular Theorem \ref{th:tgmopt}---we will show that is sufficient to choose  $\mathbf{p}$ such that
\begin{enumerate}
\item [$(i)$]
\begin{small}
\begin{equation*}  
	\mathbf{p}(\theta)^{H}\mathbf{p}(\theta)+\mathbf{p}(\theta+\pi)^{H}\mathbf{p}(\theta+\pi)>0 \quad\forall \theta\in[0,2\pi), 
\end{equation*}
\end{small}
which implies that the trigonometric function
\begin{equation}\label{eqn:def_s}
\textbf{s}(\theta) = \mathbf{p}(\theta)\left(\mathbf{p}(\theta)^{H}\mathbf{p}(\theta)+\mathbf{p}(\theta+\pi)^{H}\mathbf{p}(\theta+\pi)\right)^{-1}\mathbf{p}(\theta)^{H}
\end{equation}
is well-defined for all $\theta\in[0,2\pi)$,
\item [$(ii)$] \begin{small}
\[	\textbf{s}(\theta_0)q_{\bar{\jmath}}(\theta_0) = q_{\bar{\jmath}}(\theta_0),\]
\end{small}
\item [$(iii)$] \begin{small}
 \[
\lim_{\theta\rightarrow \theta_0} \lambda_{\bar{\jmath}}(\mathbf{f}(\theta))^{-1}(1-\lambda_{\bar{\jmath}}(\textbf{s}(\theta)))=c, \quad c\in\mathbb{R}.
\]
\end{small}
\end{enumerate}
Note that the first condition does not depend on $\textbf{f}$ and its spectral properties, then it provides a certain freedom in the choice of the grid transfer operator with respect to the problem. The second and third conditions depend on the eigenvector associated to the singularity of $\textbf{f}$, which is known, and the behaviour of the minimal eigenvalue function of  $\textbf{f}$. The latter is a scalar-valued function and its analytic properties can be investigated or approximated with the preferred mathematical tools.  
}


\subsection{TGM optimality}
\label{sec:tgmoptimality}
The following theorem proves that conditions $(i)-(iii)$ imply the approximation property $(c)$.
Combining this result with the smoothing property proved in \cite{DFFSS}, the optimality of the TGM follows from Theorem \ref{teoconv}.

\begin{theorem} \label{th:tgmopt}
%
%

	Consider the matrix $A_{N}:=\mathcal{A}_{n}(\mathbf{f})$, with $n$ even and $\mathbf{f}\in\mathcal{M}_{d}$ matrix-valued trigonometric polynomial, $\mathbf{f}\geq0$, such that condition (\ref{eqn:condition_on_f}) is satisfied. 
	Let $P^d_{n,k}$ 
	be the projecting operator defined as in equation (\ref{eqn:def_Pnkd})
	with $\mathbf{p}\in\mathcal{M}_{d}$ trigonometric polynomial satisfying conditions $(i)-(iii)$.
	Then, there exists a positive value $\gamma$ independent of $n$ such that inequality $(c)$ in Theorem \ref{teoconv} is satisfied.
\begin{proof}
The first part of the proof takes inspiration from \cite[Theorem 5.2]{DFFSS}. We report all the details for completeness, uniforming the notation.
We remind that in order to prove that there exists $\gamma>0$ independent of $n$ such that for any $x_{N}\in\mathbb{C}^{N}$
\begin{small}
	\begin{align}\label{condW}
	\min_{y\in\mathbb{C}^{K}}\|x_{N}-P^d_{n,k}y\|_{2}^{2}\leq \gamma\|x_{N}\|_{A_{N}}^{2},
	\end{align}
\end{small}
	we can choose a special instance of $y$ in such a way that the previous inequality is
	reduced to a matrix inequality in the sense of the partial ordering of the real space of
	 Hermitian matrices. For any $x_N\in\mathbb{C}^N$, let $\overline{y}\equiv\overline{y}(x_{N})\in\mathbb{C}^{K}$ be defined as
	$
	\overline{y}=[(P^d_{n,k})^{H}P^d_{n,k}]^{-1}(P^d_{n,k})^{H}x_{N}.
$
	From condition 
	$(i)$ and \cite[Proposition 4.2]{DFFSS}, it is straightforward that $(P^d_{n,k})^{H}P^d_{n,k}$ is invertible.	Therefore, (\ref{condW}) is implied by
	\begin{align*}
	\|x_{N}-P^d_{n,k}\overline{y}\|_{2}^{2}\leq \gamma\|x_{N}\|_{A_{N}}^{2},
	\end{align*}
	where the latter is equivalent to the matrix inequality $
	G_{N}(\mathbf{p})^{H}G_{N}(\mathbf{p})\leq \gamma A_{N}.$	with $G_{N}(\mathbf{p})=I_{N}-P^d_{n,k}[(P^d_{n,k})^{H}P^d_{n,k}]^{-1}(P^d_{n,k})^{H}$. By construction, the matrix $G_N(p)$ is a Hermitian unitary projector, in fact $G_N(\mathbf{p})^{H}G_N(\mathbf{p})= G_{N}(\mathbf{p})^{2} = G_N(\mathbf{p})$. As a consequence, the preceding matrix inequality can be rewritten as
\begin{small}
	\begin{align}\label{rela}
	G_{N}(\mathbf{p})\leq \gamma \mathcal{A}_{n}(\mathbf{f}).
	\end{align}
\end{small}
	We notice that $\left(K^{Odd}_{n,k}\right)^TF_{n}=\frac{1}{\sqrt{2}}F_{k}I_{n,2}$, where $I_{n,2}=\left[I_{k} | I_{k}\right]_{k\times n}$. Since we can decompose the block-circulant matrix $\mathcal{A}_{n}(\mathbf{p})=(F_n\otimes I_{d}) \diagi(\mathbf{p}(\theta_i^{(n)}))(F_{n}^{H}\otimes I_{d})$, we have
\begin{small}
	\begin{align*}
	(P^d_{n,k})^{H}=\frac{1}{\sqrt{2}}(F_{k}\otimes I_{d})(I_{n,2}\otimes I_d)\diagi(\mathbf{p}(\theta_i^{(n)})^{H})(F_{n}^{H}\otimes I_d),
	\end{align*}
\end{small}
	and the matrix $(F_{n}^{H}\otimes I_d)G_{N}(\mathbf{p})(F_{n}\otimes I_{d})$ becomes
\begin{small}
	\begin{align*}
	(F_{n}^{H}\otimes I_d)G_{N}(\mathbf{p})(F_{n}\otimes I_{d})&=I_{N}-
	\diagi(\mathbf{p}(\theta_i^{(n)}))(I_{n,2}^{T}\otimes I_d)\\	
	&\left[\diagik(\mathbf{p}(\theta_{i}^{(n)})^{H}\mathbf{p}(\theta_{i}^{(n)})+\mathbf{p}(\theta_{\tilde{\imath}}^{(n)})^{H}\mathbf{p}(\theta_{\tilde{\imath}}^{(n)}))\right]^{-1}\\
	&(I_{n,2}\otimes I_d)\diagi(\mathbf{p}(\theta_i^{(n)})^{H})
	\end{align*}
\end{small}
	where $\tilde{\imath}=i+k$. Now, it is clear that there exists a suitable permutation by rows and columns of $(F_{n}^{H}\otimes I_d)G_{N}(\mathbf{p})(F_{n}\otimes I_{d})$ such that we can obtain a $2d\times 2d$ block-diagonal matrix of the form
	\begin{small}
	\begin{equation*}
	I_{N}-\diagik\left[\begin{array}{c} 
	\mathbf{p}(\theta_{i}^{(n)}) \\ \mathbf{p}(\theta_{\tilde{\imath}}^{(n)})
	\end{array}\right]\!
	\left[\begin{array}{c} 
	(\mathbf{p}(\theta_{i}^{(n)})^{H}\mathbf{p}(\theta_{i}^{(n)})+\mathbf{p}(\theta_{\tilde{\imath}}^{(n)})^{H}\mathbf{p}(\theta_{\tilde{\imath}}^{(n)}))^{-1}
	\end{array}\right]\!
	\left[\begin{array}{cc} 
	\mathbf{p}(\theta_{i}^{(n)})^{H} & \mathbf{p}(\theta_{\tilde{\imath}}^{(n)})^{H}
	\end{array}\right].
	\end{equation*}
	\end{small}
	Therefore, by considering the same permutation by rows and columns of $(F_{n}^{H}\otimes I_d)\mathcal{A}_{n}(\mathbf{f})(F_{n}\otimes I_{d})= \diagi(\mathbf{f}(\theta_{i}^{(n)}))$, condition (\ref{rela}) is equivalent to requiring that there exists $\gamma>0$ independent of $n$ such that, $\forall j=0,\ldots,k-1$
	\begin{small}
	\begin{align*}
	&I_{2d}-\left[\begin{array}{c} 
	\mathbf{p}(\theta_{i}^{(n)}) \\ \mathbf{p}(\theta_{\tilde{\imath}}^{(n)})
	\end{array}\right]
	\left[\begin{array}{c} 
	(\mathbf{p}(\theta_{i}^{(n)})^{H}\mathbf{p}(\theta_{i}^{(n)})+\mathbf{p}(\theta_{\tilde{\imath}}^{(n)})^{H}\mathbf{p}(\theta_{\tilde{\imath}}^{(n)}))^{-1}
	\end{array}\right]
	\left[\begin{array}{cc} 
	\mathbf{p}(\theta_{i}^{(n)})^{H} & \mathbf{p}(\theta_{\tilde{\imath}}^{(n)})^{H}
	\end{array}\right] \\
	&\leq \gamma
	\left[\begin{array}{cc} 
	\mathbf{f}(\theta_{i}^{(n)}) & \\
	& \mathbf{f}(\theta_{\tilde{\imath}}^{(n)})
	\end{array}\right].
	\end{align*}
	\end{small}
	We define the set $\Omega(\theta_0)=\{\theta_0,\theta_0+\pi\}.$
	Due of the continuity of $\mathbf{p}$ and $\mathbf{f}$ it is clear that the preceding set of inequalities can be reduced to requiring that a unique inequality of the form
\begin{small}
\begin{align*}
I_{2d}-\left[\begin{array}{c} 
\mathbf{p}(\theta) \\ \mathbf{p}(\theta+\pi)
\end{array}\right]
\left[\begin{array}{c} 
(\mathbf{p}(\theta)^{H}\mathbf{p}(\theta)+\mathbf{p}(\theta+\pi)^{H}\mathbf{p}(\theta+\pi))^{-1}
\end{array}\right]\\
\left[\begin{array}{cc} 
\mathbf{p}(\theta)^{H} & \mathbf{p}(\theta+\pi)^{H}
\end{array}\right]\leq \gamma
\left[\begin{array}{cc} 
\mathbf{f}(\theta) & \\
& \mathbf{f}(\theta+\pi)
\end{array}\right]
\end{align*}
\end{small}
holds for all $\theta\in[0,2\pi)\backslash \Omega(\theta_0)$. 
Let us define $\textbf{q}(\theta)=(\mathbf{p}(\theta)^{H}\mathbf{p}(\theta)+\mathbf{p}(\theta+\pi)^{H}\mathbf{p}(\theta+\pi))^{-1}$. By simple computations, the previous inequality becomes 
\begin{small}
\begin{align}\label{eqn:after_commutativity}
I_{2d}-
\left[\begin{array}{cc} 
\mathbf{p}(\theta)\textbf{q}(\theta)\mathbf{p}(\theta)^{H} & \mathbf{p}(\theta)\textbf{q}(\theta)\mathbf{p}(\theta+\pi)^{H} \\
\mathbf{p}(\theta+\pi)\textbf{q}(\theta)\mathbf{p}(\theta)^{H} & \mathbf{p}(\theta+\pi)\textbf{q}(\theta)\mathbf{p}(\theta+\pi)^{H}
\end{array}\right]\leq \gamma
\left[\begin{array}{cc} 
\mathbf{f}(\theta) & \\
& \mathbf{f}(\theta+\pi)
\end{array}\right].
\end{align}
\end{small}
Let us define the matrix-valued function
\begin{small}
\begin{multline*}
R(\theta)=
\begin{bmatrix}
\mathbf{f}(\theta) & \\
& \mathbf{f}(\theta+\pi)
\end{bmatrix} ^{-\frac{1}{2}}
\left[\begin{array}{cc} 
I_d-\mathbf{p}(\theta)\textbf{q}(\theta)\mathbf{p}(\theta)^{H} & -\mathbf{p}(\theta)\textbf{q}(\theta)\mathbf{p}(\theta+\pi)^{H} \\
-\mathbf{p}(\theta+\pi)\textbf{q}(\theta)\mathbf{p}(\theta)^{H} & I_d-\mathbf{p}(\theta+\pi)\textbf{q}(\theta)\mathbf{p}(\theta+\pi)^{H}
\end{array}\right]\cdot\\
\cdot\begin{bmatrix}
\mathbf{f}(\theta) & \\
& \mathbf{f}(\theta+\pi)
\end{bmatrix}^{-\frac{1}{2}}.
\end{multline*}
\end{small}
Applying the Sylvester inertia law \cite{GV}, we have that the relation \eqref{eqn:after_commutativity} is verified if 
\begin{small}
\begin{align}\label{eqn:inequality_on_R}
R(\theta)\leq \gamma I_{2d} 
\end{align}
\end{small}
is satisfied. If we prove that for every $\theta\in[0,2\pi)\backslash \Omega(\theta_0)$ the matrix $R(\theta)$ is uniformly bounded in the spectral norm, then we have that there exists $\gamma>0$ which bounds the spectral radius of $R(\theta)$ and then the latter implies inequality (\ref{eqn:inequality_on_R}). To show that the matrix $R(\theta)$ is uniformly bounded in the spectral norm, we can rewrite $R(\theta)$ in components as
\begin{small}
\begin{align*}
&R(\theta)=
\begin{bmatrix}
R_{1,1}(\theta)& R_{1,2}(\theta)\\
R_{2,1}(\theta)&R_{2,2}(\theta)
\end{bmatrix}=\\
&\begin{bmatrix}
\mathbf{f}^{-\frac{1}{2}}(\theta)(I_d-\mathbf{p}(\theta)\textbf{q}(\theta)\mathbf{p}(\theta)^H)\mathbf{f}^{-\frac{1}{2}}(\theta) & \!-\mathbf{f}^{-\frac{1}{2}}(\theta)\mathbf{p}(\theta)\textbf{q}(\theta)\mathbf{p}(\theta\!+\!\pi)^H\mathbf{f}^{-\frac{1}{2}}(\theta\!+\!\pi) \\
-\mathbf{f}^{-\frac{1}{2}}(\theta\!+\!\pi)\mathbf{p}(\theta\!+\!\pi)\textbf{q}(\theta)\mathbf{p}(\theta)^H\mathbf{f}^{-\frac{1}{2}}(\theta) & \!\mathbf{f}^{-\frac{1}{2}}(\theta\!+\!\pi)(I_d\!-\!\mathbf{p}(\theta\!+\!\pi)\textbf{q}(\theta)\mathbf{p}(\theta\!+\!\pi)^H)\mathbf{f}^{-\frac{1}{2}}(\theta\!+\!\pi)
\end{bmatrix}\!\!.
\end{align*}
\end{small}
The function $\|R(\theta)\|_2:[0,2\pi)\backslash \Omega(\theta_0) \rightarrow \mathbb{R}$ is continuous and, in order to show that $R(\theta)$ is uniformly bounded in the spectral norm, Weierstrass \textcolor{red}{theorem} implies that it is sufficient to prove that the following limits exist and are finite:
\begin{small}
\begin{equation*}
\lim_{\theta\rightarrow\theta_0} \|R(\theta)\|_2, \qquad \lim_{\theta\rightarrow\theta_0+\pi} \|R(\theta)\|_2.
\end{equation*}
\end{small}
By definition, $R(\theta)$ is a Hermitian matrix for $\theta\in[0,2\pi)\backslash \Omega(\theta_0)$. Moreover, by direct computation, one can verify that the matrix on the left-hand side of (\ref{eqn:after_commutativity}) is a projector, \textcolor{red}{having} eigenvalues $0$ and $1$. Consequently, from the Sylvester inertia law, it  follows that $R(\theta)$ is a non-negative definite matrix.
We remark that in order to bound the spectral norm of a non-negative definite matrix-valued function, it is sufficient to bound its trace. Hence, we check that the spectral norms of the elements on the block diagonal of $R(\theta)$ are bounded. The latter is equivalent to verify that the  limits  
\begin{small}
\begin{equation}\label{eqn:limit_on_R11} 
\lim_{\theta\rightarrow\theta_0} \|R_{1,1}(\theta)\|_2,
\end{equation}
\begin{equation}\label{eqn:limit_on_R22} 
\quad \lim_{\theta\rightarrow\theta_0} \|R_{2,2}(\theta)\|_2,
\end{equation}
\begin{equation}\label{eqn:limit_on_R11_pi} 
\lim_{\theta\rightarrow\theta_0+\pi} \|R_{1,1}(\theta)\|_2,
\end{equation}
\begin{equation}\label{eqn:limit_on_R22_pi} 
\quad \lim_{\theta\rightarrow\theta_0+\pi} \|R_{2,2}(\theta)\|_2
\end{equation}
\end{small}
exist and are finite, which in practice requires only the proof of (\ref{eqn:limit_on_R11}). Indeed, the finiteness of (\ref{eqn:limit_on_R22}) and (\ref{eqn:limit_on_R11_pi}) is implied by the hypotheses on $\mathbf{f}$, which is non-singular in $\theta_0+\pi$. The finiteness of (\ref{eqn:limit_on_R22_pi}) can be proven as (\ref{eqn:limit_on_R11}) taking into account that $R(\theta)$ is $2\pi$-periodic.
To prove (\ref{eqn:limit_on_R11}) we note that for all $\theta\in[0,2\pi)\backslash \Omega(\theta_0)$, we can write
\begin{small}
\begin{align*}
\left\|R_{1,1}(\theta)\right\|_2=&\left\|\mathbf{f}^{-\frac{1}{2}}(\theta)(I_d-\mathbf{p}(\theta)\textbf{q}(\theta)\mathbf{p}(\theta)^H)\mathbf{f}^{-\frac{1}{2}}(\theta)\right\|_2= \left\|\mathbf{f}^{-{1}}(\theta)-\mathbf{f}^{-\frac{1}{2}}(\theta)\textbf{s}(\theta)\mathbf{f}^{-\frac{1}{2}}(\theta)\right\|_2,
\end{align*}
\end{small}
with $\textbf{s}(\theta)$ defined as in \eqref{eqn:def_s}.
Without loss of generality, we can assume that $\bar{\jmath}=1$, that is $q_1(\theta_0)$ is the eigenvector of $\mathbf{f}(\theta_0)$ associated with the eigenvalue $0$. Indeed, if $\bar{\jmath}\neq 1$, it is sufficient to permute rows and columns of $D(\theta_0)$ in the factorization in (\ref{eq:dec_f}) via a permutation matrix $\Pi$ which brings the diagonalization of $\mathbf{f}(\theta_0)$ into the desired form. Moreover, we can assume that $\|q_1(\theta_0)\|_2=1$.
From condition $(i)$ 
we have that the matrix-valued function $\textbf{s}(\theta)$ is Hermitian for all $\theta\in[0,2\pi)$. 
In addition, from condition $(ii)$
and from the latter assumption on $\bar{\jmath}$, the matrix $\textbf{s}(\theta)$ can be decomposed as $\textbf{s}(\theta)=W_s(\theta)D_s(\theta)W_s^H(\theta)$ and
\begin{small}
\begin{equation*}
\textbf{s}(\theta_0)=\left[
\begin{array}{@{\,}c@{\;}|@{\;}c@{\;}|@{\;}c@{\;}|@{\;}c@{\,}}
q_1(\theta_0)&w_2(\theta_0)&\dots & w_d(\theta_0)
\end{array}
\right]\begin{bmatrix}
1\\
& \lambda_2(\textbf{s}(\theta_0))\\
& & \ddots\\
& & & \lambda_d(\textbf{s}(\theta_0))
\end{bmatrix}
\left[
\begin{array}{@{\,}c@{\,}}
{q_1}^H(\theta_0) \\
\hline
{w_2}^H(\theta_0)\\
\hline
\vdots\\
\hline
{w_d}^H(\theta_0)
\end{array}
\right].
\end{equation*}
\end{small}
Then, we can rewrite the quantity to bound as follows: 
\begin{small}
\begin{equation*}
\begin{split}
&\lim_{\theta\rightarrow\theta_0}\|Q(\theta)D^{-1}(\theta)Q^H(\theta)-Q(\theta)D^{-\frac{1}{2}}(\theta)Q^H(\theta)
W_s(\theta)D_s(\theta)W^H_s(\theta)
Q^H(\theta)D^{-\frac{1}{2}}(\theta)Q^H(\theta)\|_2\\
&=\lim_{\theta\rightarrow\theta_0}\|D^{-1}(\theta)-D^{-\frac{1}{2}}(\theta)Q^H(\theta)
W_s(\theta)D_s(\theta)W^H_s(\theta)
Q^H(\theta)D^{-\frac{1}{2}}(\theta)\|_2.\\
\end{split}
\end{equation*}
\end{small}
By definition of $Q(\theta_0)$ and $W_s(\theta_0)$, the vector $q_{0}(\theta_0)$ is orthogonal with respect to both  $q_j(\theta_0)$, $w_j(\theta_0)$, $j=2,\dots,d$. Denoting by $\textbf{0}_{d-1}$ the null row vector of size $d-1$, we have
\begin{small}
\begin{equation*}
\lim_{\theta\rightarrow\theta_0}Q^H(\theta)
W_s(\theta)=\left[
\begin{array} { c c }
                \begin{array}{c}
                q_1(\theta_0)^Hq_1(\theta_0)
                \end{array}                   & \textbf{0}_{d-1}\\
                \textbf{0}_{d-1}^T                             & \begin{array}{cc}
                  M(\theta_0)
                                               \end{array}
\end{array}
\right],
\end{equation*}
\end{small}
where $M(\theta)$ is a matrix-valued function which is well-defined and continuous on $[0,2\pi]$.
Then, since the eigenvalue functions $\lambda_i(\mathbf{f}(\theta))^{-1}$, for $i=2,\dots,d$, are well-defined and continuous on $[0,2\pi],$ see Lemma \ref{lemm:continuity_eig}, the quantity to bound becomes
\begin{small}
\begin{equation*}
\begin{split}
&\left\|\!  \begin{bmatrix}\!
          \underset{\theta\rightarrow\theta_0}{\lim}\lambda_1(\mathbf{f}(\theta))^{-1}(1-\lambda_1(\textbf{s}(\theta)))
               & \textbf{0}_{d-1}\\
               \textbf{0}_{d-1}^T                             &
               \begin{bmatrix}
\lambda_2(\mathbf{f}(\theta_0))^{-1}\\
 & \!\!\!\!\ddots\!\!\\
 & &\lambda_d(\mathbf{f}(\theta_0))^{-1}
\end{bmatrix}\!\left(I_{d-1}\!-\!M(\theta_0)M^T(\theta_0)\right)\!
\end{bmatrix}\!\right\|_2
\end{split}\!.
\end{equation*}
\end{small}
Consequently, the thesis follows from condition $(iii)$. 
\end{proof}
\end{theorem}
\begin{remark}
The proof of Theorem \ref{th:tgmopt} requires that $\theta_0 \neq \theta_i^{(n)}$ for all $i$ and $n$.
Nevertheless, in practice, our multigrid method works well even if $\theta_0 = \theta_i^{(n)}$ for a certain $i$ and $n$, as it can happen in the circulant case. In such case the coefficient matrix $A_N$ is singular but the multigrid method converges anyway, since it works on the orthogonal complement of the eigenvector corresponding to the zero eigenvalue.
We could add a rank one correction, like the one used for scalar symbols in \cite{ADS}, but this would only lead to unnecessary complication of notation, see \cite{AD}.
\end{remark}
{
In practical applications choosing a $\mathbf{p}$ such that conditions $(ii)$ and $(iii)$ are verified could not be  trivial. 
Hence, in the following, assuming that $\mathbf{p}$ satisfies the condition $(i)$ 
so that the matrix-valued function $\textbf{s}$ is well-defined,
we provide two useful results, Lemma \ref{lemm:condition_on_s}-\ref{lemm:condition_on_s_p_invertible}, that can be used to construct $\mathbf{p}$ that fulfills condition $(ii)$. 
Analogously, Lemma \ref{lemm:condition_iii} shows how to deal with condition $(iii)$ under some additional hypotheses on $\mathbf{p}$ and $\mathbf{f}$.}

\begin{lemma}\label{lemm:condition_on_s}
Let $\mathbf{f}$ be a matrix-valued trigonometric polynomial, $\mathbf{f}\geq0$ that satisfies condition (\ref{eqn:condition_on_f}). 
Assume $\mathbf{p}$ is a matrix-valued trigonometric polynomial such that condition $(i)$ is fulfilled, so that the matrix-valued function $\textbf{s}$ defined as in \eqref{eqn:def_s} is well-defined. 
Assume that the eigenvector $q_{\bar{\jmath}}(\theta_0)$ associate with the ill-conditioned subspace of $\mathbf{f}(\theta_0)$,
i.e., $\mathbf{f}(\theta_0)q_{\bar{\jmath}}(\theta_0)=0 q_{\bar{\jmath}}(\theta_0)$, is such that:
\begin{enumerate}
\item $q_{\bar{\jmath}}(\theta_0)$ is an eigenvector of $\mathbf{p}(\theta_0)$, associated to $\lambda^{(1)}_{\bar{\jmath}}\neq 0$ that is
$$\mathbf{p}(\theta_0)q_{\bar{\jmath}}(\theta_0)=\lambda^{(1)}_{\bar{\jmath}}q_{\bar{\jmath}}(\theta_0);$$
\item $q_{\bar{\jmath}}(\theta_0)$ is an eigenvector of $\mathbf{p}(\theta_0+\pi)$ associated with the zero eigenvalue, that is
$$\mathbf{p}(\theta_0+\pi)q_{\bar{\jmath}}(\theta_0)=0 q_{\bar{\jmath}}(\theta_0);$$
\item $q_{\bar{\jmath}}(\theta_0)$ is an eigenvector of $\mathbf{p}(\theta_0)^H$, associated to $\lambda^{(2)}_{\bar{\jmath}}\neq 0$, that is $$\mathbf{p}(\theta_0)^Hq_{\bar{\jmath}}(\theta_0)=\lambda^{(2)}_{\bar{\jmath}}q_{\bar{\jmath}}(\theta_0).$$
\end{enumerate}
Then condition $(ii)$ is satisfied.
\begin{proof}
From all the hypotheses on $q_{\bar{\jmath}}(\theta_0)$ and by direct computation, we have
$
\left(\mathbf{p}(\theta_0)^{H}\mathbf{p}(\theta_0)+\mathbf{p}(\theta_0+\pi)^H\mathbf{p}(\theta_0+\pi)\right)q_{\bar{\jmath}}(\theta_0)=
\lambda^{(1)}_{\bar{\jmath}}\lambda^{(2)}_{\bar{\jmath}} q_{\bar{\jmath}}(\theta_0).
$
Then, by definition of $\textbf{s}(\theta)$ in \eqref{eqn:def_s}, it holds that
\begin{small}
\begin{align*}
\textbf{s}(\theta_0)q_{\bar{\jmath}}(\theta_0)&=\mathbf{p}(\theta_0)\left(\mathbf{p}(\theta_0)^{H}\mathbf{p}(\theta_0)+\mathbf{p}(\theta_0+\pi)^H\mathbf{p}(\theta_0+\pi)\right)^{-1}\mathbf{p}(\theta_0)^Hq_{\bar{\jmath}}(\theta_0)\\
&=\lambda_{\bar{\jmath}}^{(2)}\mathbf{p}(\theta_0)\left(\mathbf{p}(\theta_0)^{H}\mathbf{p}(\theta_0)+\mathbf{p}(\theta_0+\pi)^H\mathbf{p}(\theta_0+\pi)\right)^{-1}q_{\bar{\jmath}}(\theta_0)\\
&=\lambda_{\bar{\jmath}}^{(2)}\frac{1}{\lambda^{(1)}_{\bar{\jmath}}\lambda^{(2)}_{\bar{\jmath}}}\mathbf{p}(\theta_0)q_{\bar{\jmath}}(\theta_0)
=q_{\bar{\jmath}}(\theta_0).
\end{align*}
\end{small}
\end{proof}
\end{lemma}
The next lemma provides other sufficient conditions to verify hypothesis $(ii)$, if the projector is associated with a trigonometric polynomial $\mathbf{p}$ which is non-singular in the considered point $\theta_0$. 
\begin{lemma}\label{lemm:condition_on_s_p_invertible}
With the assumption and notation of Lemma \ref{lemm:condition_on_s} where the condition 3. is replaced with 
\begin{enumerate}
\item [3 bis.] $\mathbf{p}(\theta_0)$ is non-singular.
\end{enumerate}
Then condition $(ii)$ is satisfied.
\begin{proof}
By definition of $\textbf{s}$ in equation \eqref{eqn:def_s}, we have
\begin{small}
\begin{equation*}
\begin{split}
\textbf{s}(\theta_0)^{-1}&=\left(\mathbf{p}(\theta_0)\left(\mathbf{p}(\theta_0)^{H}\mathbf{p}(\theta_0)+\mathbf{p}(\theta_0+\pi)^H\mathbf{p}(\theta_0+\pi)\right)^{-1}\mathbf{p}(\theta_0)^H\right)^{-1}\\
&=\mathbf{p}(\theta_0)^{-H}\left(\mathbf{p}(\theta_0)^{H}\mathbf{p}(\theta_0)+\mathbf{p}(\theta_0+\pi)^H\mathbf{p}(\theta_0+\pi)\right)^{-1}\mathbf{p}(\theta_0)^{-1}\\
&= I_d+\mathbf{p}(\theta_0)^{-H}\mathbf{p}(\theta_0+\pi)^H\mathbf{p}(\theta_0+\pi)\mathbf{p}(\theta_0)^{-1}.
\end{split}
\end{equation*}
\end{small}
Then, it holds
\begin{small}
\begin{equation*}
\begin{split}
s(\theta_0)^{-1}q_{\bar{\jmath}}(\theta_0)=&\left[I_d+\mathbf{p}(\theta_0)^{-H}\mathbf{p}(\theta_0+\pi)^H\mathbf{p}(\theta_0+\pi)\mathbf{p}(\theta_0)^{-1}\right]q_{\bar{\jmath}}(\theta_0)=\\
&q_{\bar{\jmath}}(\theta_0)+\frac{1}{\lambda_{\bar{\jmath}}^{(1)}}\mathbf{p}(\theta_0)^{-H}\mathbf{p}(\theta_0+\pi)^H\mathbf{p}(\theta_0+\pi)q_{\bar{\jmath}}(\theta_0)=q_{\bar{\jmath}}(\theta_0)
\end{split}
\end{equation*}
\end{small}
and then condition $(ii)$ follows.
\end{proof}
\end{lemma}
 {
\textcolor{red}{Finally, we present Lemma \ref{lemm:condition_iii} to  simplify the validation of condition $(iii)$  in possible applications. This provides some additional hypotheses on \textcolor{red}{$\textbf{p}$} and \textcolor{red}{$\textbf{f}$} that can be considered when proving that condition $(iii)$ is satisfied.  For the aforementioned purpose we first introduce the following remark containing algebraic calculations in order to speed up the proof of the lemma. 
}
\begin{remark}\label{rmk:eig_p_pi}
Suppose that we can write $\mathbf{p}(\theta)$ as
\begin{small}
\begin{equation*}
\mathbf{p}(\theta)=\left[
\begin{array}{@{\,}c@{\;}|@{\;}c@{\;}|@{\;}c@{\;}|@{\;}c@{\,}}
h_1(\theta)&h_2(\theta)&\dots & h_d(\theta)
\end{array}
\right]\begin{bmatrix}
\lambda_{1}(\mathbf{p}(\theta)) & \textbf{0}_{d-1}\\
\textbf{0}_{d-1}^T& M(\theta)
\end{bmatrix}
\left[
\begin{array}{@{\,}c@{\,}}
{h_1}^H(\theta) \\
\hline
{h_2}^H(\theta)\\
\hline
\vdots\\
\hline
{h_d}^H(\theta)
\end{array}
\right],
\end{equation*}
\end{small} 
where $M(\theta)\in \mathbb{C}^{(d-1)\times (d-1)}, 
$  for each $\theta$,  $\lim_{\theta\rightarrow\theta_0}h_1(\theta+\pi)=q_{\bar{\jmath}}(\theta_0)$, and $h_{\jmath}(\theta_0+\pi)^H q_{\bar{\jmath}}(\theta_0)=0,$ for ${\jmath}=2,\dots, d$. Then, we can write
\begin{small}
\begin{equation*}
\begin{split}
&\lim_{\theta\rightarrow\theta_0} \mathbf{p}(\theta+\pi)q_{\bar{\jmath}}(\theta_0)=\\
&\lim_{\theta\rightarrow \theta_0}\!\left[
\begin{array}{@{\,}c@{\;}|@{\;}c@{\;}|@{\;}c@{\;}|@{\;}c@{\,}}
\!h_1(\theta+\pi)&h_2(\theta+\pi)&\dots & h_d(\theta+\pi)\!
\end{array}
\right]\!\!\begin{bmatrix}
\lambda_{1}(\mathbf{p}(\theta+\pi)) & \textbf{0}_{d-1}\\
\textbf{0}_{d-1}^T& \!\!M(\theta+\pi)
\end{bmatrix}
\!\!\left[
\begin{array}{@{\,}c@{\,}}
{h_1}^H(\theta+\pi) \\
\hline
{h_2}^H(\theta+\pi)\!\\
\hline
\vdots\\
\hline
{h_d}^H(\theta+\pi)\end{array}
\right]\!\!q_{\bar{\jmath}}(\theta_0)\!
=\\
& \lim_{\theta\rightarrow \theta_0}\! \left[
\begin{array}{@{\,}c@{\;}|@{\;}c@{\;}|@{\;}c@{\;}|@{\;}c@{\,}}
\!h_1(\theta+\pi)&h_2(\theta+\pi)&\dots & h_d(\theta+\pi)\!
\end{array}
\right]\!\!\begin{bmatrix}
\lambda_{1}(\mathbf{p}(\theta+\pi)) & \textbf{0}_{d-1}\\
\textbf{0}_{d-1}^T& \!\!M(\theta+\pi)
\end{bmatrix}
\!\!\left[
\begin{array}{@{\,}c@{\,}}
{q_{\bar{\jmath}}}^H(\theta_0) \\
\hline
{h_2}^H(\theta_0+\pi)\\
\hline
\vdots\\
\hline
{h_d}^H(\theta_0+\pi)\end{array}
\right]\!\!q_{\bar{\jmath}}(\theta_0)\!=\\
&  \lim_{\theta\rightarrow \theta_0} \left[
\begin{array}{@{\,}c@{\;}|@{\;}c@{\;}|@{\;}c@{\;}|@{\;}c@{\,}}
\!h_1(\theta+\pi)&h_2(\theta+\pi)&\dots & h_d(\theta+\pi)\!
\end{array}
\right]\!\!\begin{bmatrix}
\lambda_{1}(\mathbf{p}(\theta+\pi)) & \textbf{0}_{d-1}\\
\textbf{0}_{d-1}^T& \!\!M(\theta+\pi)
\end{bmatrix}\!\!\begin{bmatrix}
1\\
0\\
\vdots\\
0
\end{bmatrix}=\\
&  \lim_{\theta\rightarrow \theta_0} \left[
\begin{array}{@{\,}c@{\;}|@{\;}c@{\;}|@{\;}c@{\;}|@{\;}c@{\,}}
h_1(\theta+\pi)&h_2(\theta+\pi)&\dots & h_d(\theta+\pi)
\end{array}
\right] \begin{bmatrix}
\lambda_{1}(\mathbf{p}(\theta+\pi))\\
0\\
\vdots\\
0
\end{bmatrix}=\\
&  \lim_{\theta\rightarrow \theta_0} \lambda_{1}(\mathbf{p}(\theta+\pi)) h_1(\theta+\pi)=\lim_{\theta\rightarrow \theta_0} \lambda_{1}(\mathbf{p}(\theta+\pi)) q_{\bar{\jmath}}(\theta_0).
\end{split}
\end{equation*}
\end{small}
\end{remark}

\begin{lemma}\label{lemm:condition_iii}
Assume that $\mathbf{p}$ and $\mathbf{f}$ are matrix-valued functions which satisfy the requirements of the Lemma~\ref{lemm:condition_on_s}. If 
\begin{enumerate}
\item {\small $$\lambda_{\bar{\jmath}}^{(1)}=\lambda_{\bar{\jmath}}^{(2)}=\lambda_{\bar{\jmath}},$$}
\item \begin{small}
$$
\lim_{\theta\rightarrow 	\theta_0} \frac{|\lambda_{\bar{\jmath}}(\mathbf{p}(\theta+\pi)|^2}{\lambda_{\bar{\jmath}}(\mathbf{f}(\theta))}=c,
$$
\end{small}
then, the condition $(iii)$ of the Theorem \ref{th:tgmopt} is satisfied. That is,
{\small $$ \lim_{\theta\rightarrow 	\theta_0} \frac{1-\lambda_{\bar{\jmath}}(\textbf{s}(\theta))}{\lambda_{\bar{\jmath}}(\mathbf{f}(\theta))}=c. $$}
\end{enumerate}
\begin{proof}
From the hypotheses on $\textbf{s}(\theta)$, we have that ${\textbf{s}}(\theta)$ is equal to
\begin{small}
\begin{equation*}
\begin{split}
& \left[
\begin{array}{@{\;}c@{\;}|@{\;}c@{\;}|@{\;}c@{\;}|@{\;}c@{\;}|@{\;}c@{\;}}
w_1(\theta)&\dots &w_{\bar{\jmath}}(\theta)&\dots& w_d(\theta)
\end{array}
\right]\begin{bmatrix}
\lambda_1({\textbf{s}}(\theta))& & &  &\\
& \ddots & &  &\\
& & \lambda_{\bar{\jmath}}({\textbf{s}}(\theta))& & \\
& & & \ddots & \\
& & & &\lambda_d({\textbf{s}}(\theta))
\end{bmatrix}
\left[
\begin{array}{ccccccc}
{w_1}^H(\theta) \\
\hline
\vdots\\
\hline
w_{\bar{\jmath}}^H(\theta)\\
\hline
\vdots\\
\hline
{w_d}^H(\theta)
\end{array}
\right]
\end{split},
\end{equation*}
\end{small}
where, $\lim_{\theta\rightarrow\theta_0}w_{\jmath}(\theta)=q_{\bar{\jmath}}(\theta_0).$
By definition of $\textbf{s}$, we can write
\begin{small}
\begin{equation*}
\begin{split}
 &\lim_{\theta\rightarrow 	\theta_0} \frac{1-\lambda_{\bar{\jmath}}(\textbf{s}(\theta))}{\lambda_{\bar{\jmath}}(\mathbf{f}(\theta))}=\\
 &
 \lim_{\theta\rightarrow 	\theta_0}\frac{1-w_{\jmath}(\theta)^H\mathbf{p}(\theta)\left(\mathbf{p}(\theta)^H\mathbf{p}(\theta)+\mathbf{p}(\theta+\pi)^H\mathbf{p}(\theta+\pi)\right)^{-1}\mathbf{p}(\theta)^Hw_{\jmath}(\theta)}{\lambda_{\bar{\jmath}}(\mathbf{f}(\theta))}.
\end{split}
\end{equation*}
\end{small}
Since $\mathbf{p}(\theta_0)q_{\bar{\jmath}}(\theta_0)=\mathbf{p}(\theta_0)^Hq_{\bar{\jmath}}(\theta_0)=\lambda_{\bar{\jmath}}q_{\bar{\jmath}}(\theta_0)$,
\begin{small}
\begin{equation*}
\begin{split}
&\lim_{\theta\rightarrow 	\theta_0} \frac{1-\lambda_{\bar{\jmath}}(\textbf{s}(\theta))}{\lambda_{\bar{\jmath}}(\mathbf{f}(\theta))}=\\
& \lim_{\theta\rightarrow 	\theta_0}\frac{1-|\lambda_{\bar{\jmath}}|q_{\bar{\jmath}}(\theta_0)^H\left(\mathbf{p}(\theta)^H\mathbf{p}(\theta)+\mathbf{p}(\theta+\pi)^H\mathbf{p}(\theta+\pi)\right)^{-1}q_{\bar{\jmath}}(\theta_0)}{\lambda_{\bar{\jmath}}(\mathbf{f}(\theta))}.
\end{split}
\end{equation*}
\end{small}
Note that,
\begin{small}
\begin{equation*}
\begin{split}
&\lim_{\theta\rightarrow 	\theta_0} q_{\bar{\jmath}}(\theta_0)^H\left(\mathbf{p}(\theta)^H\mathbf{p}(\theta)+\mathbf{p}(\theta+\pi)^H\mathbf{p}(\theta+\pi)\right)q_{\bar{\jmath}}(\theta_0)=\\
&
\lim_{\theta\rightarrow 	\theta_0} q_{\bar{\jmath}}(\theta_0)^H\mathbf{p}(\theta)^H\mathbf{p}(\theta)q_{\bar{\jmath}}(\theta_0) +\lim_{\theta\rightarrow 	\theta_0}q_{\bar{\jmath}}(\theta_0)^H \mathbf{p}(\theta+\pi)^H\mathbf{p}(\theta+\pi)q_{\bar{\jmath}}(\theta_0)= \\
& \lim_{\theta\rightarrow 	\theta_0} |\lambda_{\bar{\jmath}}|^2+ q_{\bar{\jmath}}(\theta_0)^H \mathbf{p}(\theta+\pi)^H\mathbf{p}(\theta+\pi)q_{\bar{\jmath}}(\theta_0).
\end{split}
\end{equation*}
\end{small}
Then, 
\begin{small}
\begin{equation*}
\begin{split}
&\lim_{\theta\rightarrow 	\theta_0} q_{\bar{\jmath}}(\theta_0)^H\left(\mathbf{p}(\theta)^H\mathbf{p}(\theta)+\mathbf{p}(\theta+\pi)^H\mathbf{p}(\theta+\pi)\right)^{-1}q_{\bar{\jmath}}(\theta_0)=\\
&\lim_{\theta\rightarrow 	\theta_0} \frac{1}{|\lambda_{\bar{\jmath}}|^2+ q_{\bar{\jmath}}(\theta_0)^H \mathbf{p}(\theta+\pi)^H\mathbf{p}(\theta+\pi)q_{\bar{\jmath}}(\theta_0)}.
\end{split}
\end{equation*}
\end{small}
Consequently, we can write
\begin{small}
\begin{equation*}
\begin{split}
 &\lim_{\theta\rightarrow 	\theta_0} \frac{1-\lambda_{\bar{\jmath}}(\textbf{s}(\theta))}{\lambda_{\bar{\jmath}}(\mathbf{f}(\theta))}= \frac{1-\frac{|\lambda_{\bar{\jmath}}|^2}{|\lambda_{\bar{\jmath}}|^2+ q_{\bar{\jmath}}(\theta_0)^H \mathbf{p}(\theta+\pi)^H\mathbf{p}(\theta+\pi)q_{\bar{\jmath}}(\theta_0)}}{\lambda_{\bar{\jmath}}(\mathbf{f}(\theta))}=\\
 &\lim_{\theta\rightarrow 	\theta_0} \frac{q_{\bar{\jmath}}(\theta_0)^H \mathbf{p}(\theta+\pi)^H\mathbf{p}(\theta+\pi)q_{\bar{\jmath}}(\theta_0)}{\lambda_{\bar{\jmath}}(\mathbf{f}(\theta))\left(|\lambda_{\bar{\jmath}}|^2+ q_{\bar{\jmath}}(\theta_0)^H \mathbf{p}(\theta+\pi)^H\mathbf{p}(\theta+\pi)q_{\bar{\jmath}}(\theta_0)\right)}=\\
&c \lim_{\theta\rightarrow 	\theta_0} \frac{q_{\bar{\jmath}}(\theta_0)^H \mathbf{p}(\theta+\pi)^H\mathbf{p}(\theta+\pi)q_{\bar{\jmath}}(\theta_0)}{\lambda_{\bar{\jmath}}(\mathbf{f}(\theta))},
 \end{split}
\end{equation*}
\end{small}
where in the latter equality we used the fact that  $|\lambda_{\bar{\jmath}}|^2>0.$ Hence, from Remark \ref{rmk:eig_p_pi}, we have that
\begin{small}
\begin{equation*}
\lim_{\theta\rightarrow 	\theta_0} \frac{1-\lambda_{\bar{\jmath}}(\textbf{s}(\theta))}{\lambda_{\bar{\jmath}}(\mathbf{f}(\theta))}=c\lim_{\theta\rightarrow 	\theta_0} \frac{|\lambda_{\bar{\jmath}}(\mathbf{p}(\theta+\pi)|^2}{\lambda_{\bar{\jmath}}(\mathbf{f}(\theta))}. 
\end{equation*} 
\end{small}
Then, the thesis follows from hypothesis 2.  
\end{proof}
\end{lemma}
\textcolor{red}{
In conclusion, in order to prove that a considered symbol $\mathbf{p}$ is a good choice to construct the grid transfer operator for the TGM method it is sufficient to check if $\mathbf{p}$ satisfies:
\begin{enumerate}
\item condition $(i)$,
\item Lemma \ref{lemm:condition_on_s} or Lemma \ref{lemm:condition_on_s_p_invertible},
\item Lemma \ref{lemm:condition_iii}.
\end{enumerate}
}

\subsection{V-cycle optimality}\label{sec:condition_Vcycle}
{Following the proof of the Theorem \ref{th:tgmopt} and the results in \cite{NN}, it is possible to derive also conditions for the convergence and optimality of the V-cycle in the block case. 
Indeed, according to Lemma \ref{lem:V}, it is sufficient to prove that 
there exists a $\delta$ independent from $n$ such that, for all $\ell$, $
\|\pi_{A_{n_\ell}}\|_2\le \delta,$
or, equivalently, that 
$
\|A_{n_\ell}^{1/2} P_{n_\ell,k_\ell}(P_{n_\ell,k_\ell}^H A_{n_\ell}P_{n_\ell,k_\ell})^{-1}P_{n_\ell,k_\ell}A_{n_\ell}^{1/2}\|_2\le \delta,
$
which is implied by
\begin{small}
\begin{equation}\label{eq:con_notay_block}
A_{n_\ell}^{1/2} P_{n_\ell,k_\ell}(P_{n_\ell,k_\ell}^H A_{n_\ell}P_{n_\ell,k_\ell})^{-1}P_{n_\ell,k_\ell}A_{n_\ell}^{1/2} \le \delta I_N.
\end{equation}
\end{small}
}
{As proved in \cite[Proposition 1]{DFFSS}, all matrices $A_{n_\ell}$ have a block-circulant structure and share the same spectral properties, in particular for $\theta_0=0$, see also equation \eqref{f2t}. Therefore, the following analysis performed \textcolor{red}{at} the first level could be repeated unchanged at a generic level $\ell$.}
Following the steps of the proof of Theorem \ref{th:tgmopt}, condition (\ref{eq:con_notay_block}) becomes
\begin{small}
\begin{align*}
	&\diagik \left[\begin{array}{cc} 
	\mathbf{f}(\theta_{i}^{(n)}) & \\
	& \mathbf{f}(\theta_{\tilde{\imath}}^{(n)})
	\end{array}\right]^{\frac{1}{2}}\left[\begin{array}{c} 
	\mathbf{p}(\theta_{i}^{(n)}) \\ \mathbf{p}(\theta_{\tilde{\imath}}^{(n)})
	\end{array}\right]\\
	&\cdot\diagik \left[\begin{array}{c} 
	\mathbf{p}(\theta_{i}^{(n)})^{H}\mathbf{f}(\theta_{i}^{(n)})\mathbf{p}(\theta_{i}^{(n)})+\mathbf{p}(\theta_{\tilde{\imath}}^{(n)})^{H}\mathbf{f}(\theta_{\tilde{\imath}}^{(n)})\mathbf{p}(\theta_{\tilde{\imath}}^{(n)})\end{array}\right]^{-1}	\\
	&\cdot\diagik 
	\left[\begin{array}{cc} 
	\mathbf{p}(\theta_{i}^{(n)})^{H} & \mathbf{p}(\theta_{\tilde{\imath}}^{(n)})^{H}
	\end{array}\right] \left[\begin{array}{cc} 
	\mathbf{f}(\theta_{i}^{(n)}) & \\
	& \mathbf{f}(\theta_{\tilde{\imath}}^{(n)})
	\end{array}\right]^{\frac{1}{2}}\le \delta I_N.
	\end{align*}
\end{small}
The latter is equivalent to require that  are bounded the components of the matrix-valued function 
\begin{small}
\begin{equation*}
\begin{split}
&Z(\theta)=\\
&\left[\begin{array}{cc} 
\mathbf{f}(\theta) & \\
& \mathbf{f}(\theta+\pi)
\end{array}\right]^{\frac{1}{2}}\left[\begin{array}{c} 
\mathbf{p}(\theta) \\ \mathbf{p}(\theta+\pi)
\end{array}\right]
\hat{\mathbf{f}}^{-1}(2\theta)\left[\begin{array}{cc} 
\mathbf{p}(\theta)^{H} & \mathbf{p}(\theta+\pi)^{H}
\end{array}\right]
\left[\begin{array}{cc} 
\mathbf{f}(\theta) & \\
& \mathbf{f}(\theta+\pi)
\end{array}\right]^{\frac{1}{2}},
\end{split}
\end{equation*}
\end{small}
where 
\begin{small}
\begin{align}\label{f2t}
  \hat{\mathbf{f}}(\theta)=\frac{1}{2}\left(\mathbf{p}\left(\frac{\theta}{2}\right)^{H}\mathbf{f}\left(\frac{\theta}{2}\right)\mathbf{p}\left(\frac{\theta}{2}\right)+
	\mathbf{p}\left(\frac{\theta}{2}+\pi\right)^{H}\mathbf{f}\left(\frac{\theta}{2}+\pi\right)\mathbf{p}\left(\frac{\theta}{2}+\pi\right)\right)
\end{align}
\end{small}
is  the generating function of $\mathcal{A}_{k}(\hat{\mathbf{f}})=P_{n,k}^H\mathcal{A}_{n}(\mathbf{f})P_{n,k}$.
We have 
\begin{small}
\begin{equation*}
\begin{split}
 \|Z_{1,1}(\theta)\|&=\|\mathbf{f}(\theta) \mathbf{p}(\theta)\hat{\mathbf{f}}^{-1}(2\theta) \mathbf{p}(\theta)^{H}\mathbf{f}(\theta)\|\le\\
 & \le \|\mathbf{f}(\theta) \mathbf{p}(\theta)\hat{\mathbf{f}}^{-1}(2\theta)\|\|\mathbf{p}(\theta)^{H}\mathbf{f}(\theta)\|,\\
\|Z_{2,2}(\theta)\|& =\|\mathbf{f}(\theta+\pi) \mathbf{p}(\theta+\pi)\hat{\mathbf{f}}^{-1}(2\theta) \mathbf{p}(\theta)^{H}\mathbf{f}(\theta+\pi)\|\le\\
& \le \|\mathbf{f}(\theta+\pi)\|\| \mathbf{p}(\theta+\pi)\hat{\mathbf{f}}^{-1}(2\theta)\|\|\mathbf{p}(\theta+\pi)^{H}\mathbf{f}(\theta+\pi)\|,\\
 \|Z_{1,2}(\theta)\|, \|Z_{2,1}(\theta)\| &\le   \|\mathbf{f}(\theta) \mathbf{p}(\theta)\hat{\mathbf{f}}^{-1}(2\theta)\|\|\mathbf{p}(\theta+\pi)^{H}\mathbf{f}(\theta+\pi)\|.
\end{split}
\end{equation*}
\end{small}

Since $\mathbf{f}(\theta)$ and $ \mathbf{p}(\theta)$ are trigonometric polynomials, the quantities $\| \mathbf{p}(\theta)^H\mathbf{f}(\theta)\|$, $\|\mathbf{f}(\theta+\pi)\|$ and $\|\mathbf{p}(\theta+\pi)^H\mathbf{f}(\theta+\pi) \|$ are bounded. Hence,  we have to prove that
\begin{small}
\begin{equation}\label{eq:bound_1_2_vcycle}
\|\mathbf{f}(\theta) \mathbf{p}(\theta)\hat{\mathbf{f}}^{-1}(2\theta)\|<\infty,\quad
\| \mathbf{p}(\theta+\pi)\hat{\mathbf{f}}^{-1}(2\theta)\|<\infty.
\end{equation}
\end{small}
Consequently, a key point is to investigate the properties of the generating function at the coarse levels. The following lemma will be useful tools for the aforementioned purpose.

\begin{lemma}\label{lemma:properties_f_hat_f}
Let $\mathbf{f}$ be defined as in Theorem \ref{th:tgmopt} and $\hat{\mathbf{f}}$ be defined as in formula (\ref{f2t}). Assume that $q_{\bar{\jmath}}(\theta_0)$ is the eigenvector associated with the ill-conditioned subspace of $\mathbf{f}(\theta_0)$. In addition, assume that the eigenvector $q_{\bar{\jmath}}(\theta_0)$ is such that:
\begin{enumerate}
\item [(a)] $q_{\bar{\jmath}}(\theta_0)$ is an eigenvector of $\mathbf{p}(\theta_0)$, associated to $\lambda^{(1)}\neq 0$ that is $\mathbf{p}(\theta_0)q_{\bar{\jmath}}(\theta_0)=\lambda^{(1)}q_{\bar{\jmath}}(\theta_0);$
\item [(b)] $q_{\bar{\jmath}}(\theta_0)$ is an eigenvector of $\mathbf{p}(\theta_0+\pi)$ associated with the zero eigenvalue, that is $\mathbf{p}(\theta_0+\pi)q_{\bar{\jmath}}(\theta_0)=0 q_{\bar{\jmath}}(\theta_0).$
\end{enumerate}
 Then, the following properties are fulfilled:
\begin{enumerate}
\item $\hat{\mathbf{f}}$ is an Hermitian matrix-valued trigonometric polynomial;
\item $\hat{\mathbf{f}}(\theta)\ge 0$, $\forall \, \theta\in [0,2\pi]$;
\item  $\hat{\mathbf{f}}(2\theta_0)q_{\bar{\jmath}}(\theta_0)=0q_{\bar{\jmath}}(\theta_0);$
\item $\hat{\mathbf{f}}(\theta)>0,$ $\,\theta\neq 2\theta_0$;
\item  If \begin{small}
\begin{equation}\label{eq:condition_for_order_zero}
\lim_{\theta\rightarrow\theta_0} \frac{\lambda_{\bar{\jmath}}^2(\mathbf{p}(\theta +\pi))}{\lambda_{\bar{\jmath}}(\mathbf{f}(\theta))}<\infty,
\end{equation}
\end{small}
then 
{\small $$\lim_{\theta\rightarrow\theta_0}\frac{\lambda_{\bar{\jmath}}(\hat{\mathbf{f}}(2\theta))}{\lambda_{\bar{\jmath}}(\mathbf{f}(\theta))}=c, \, c\neq 0\in \mathbb{R}.$$}
\end{enumerate}
\begin{proof} $ $
\begin{enumerate}
\item It is straightforward to see that $\hat{\mathbf{f}}$ is an Hermitian matrix-valued trigonometric polynomial from its definition in (\ref{f2t}). In particular, $\hat{\mathbf{f}}$ is obtained by sums and products of the trigonometric polynomials $\mathbf{p}$ and $\mathbf{f}$.  
\item Assume $y\in \mathbb{R}^d$, then we have
\begin{small}
\[y^T\hat{\mathbf{f}}y= \frac{1}{2}y^T\mathbf{p}\left(\frac{\theta}{2}\right)^{H}\mathbf{f}\left(\frac{\theta}{2}\right)\mathbf{p}\left(\frac{\theta}{2}\right) y+\frac{1}{2}y^T
	\mathbf{p}\left(\frac{\theta}{2}+\pi\right)^{H}\mathbf{f}\left(\frac{\theta}{2}+\pi\right)\mathbf{p}\left(\frac{\theta}{2}+\pi\right)y.\] \end{small}
	By hypotheses, $\mathbf{f}(\frac{\theta}{2})$ and $\mathbf{f}(\frac{\theta}{2}+\pi)$ are non negative matrices, then $y^T\hat{\mathbf{f}}y\ge 0.$ 
	\item By definition, it holds
	$$\hat{\mathbf{f}}(2\theta_0)q_{\bar{\jmath}}=\frac{1}{2}\mathbf{p}^{H}(\theta_0)\mathbf{f}({\theta_0})\mathbf{p}(\theta_0)q_{\bar{\jmath}}+\frac{1}{2}
	\mathbf{p}^{H}(\theta_0+\pi)\mathbf{f}(\theta_0+\pi)\mathbf{p}(\theta_0+\pi)q_{\bar{\jmath}}.$$
The vector $q_{\bar{\jmath}}(\theta_0)$ is the eigenvector of $ \mathbf{p}(\theta_0+\pi)$ and $\mathbf{f}(\theta_0)$ associated with the zero eigenvalue and it is the eigenvector  of $\mathbf{p}(\theta_0)$  associated with $\lambda^{(1)}$. Then,  we have
	 $$\hat{\mathbf{f}}(2\theta_0)q_{\bar{\jmath}}=\frac{\lambda^{(1)}}{2}\mathbf{p}^{H}(\theta_0)\mathbf{f}({\theta_0})q_{\bar{\jmath}}(\theta_0)+0q_{\bar{\jmath}}(\theta_0)=0q_{\bar{\jmath}}(\theta_0).$$
	 \item	Assume $y\neq 0\in \mathbb{R}^d$, then we have
	 \begin{small}
\[y^T\hat{\mathbf{f}}y= \frac{1}{2}y^T\mathbf{p}\left(\frac{\theta}{2}\right)^{H}\mathbf{f}\left(\frac{\theta}{2}\right)\mathbf{p}\left(\frac{\theta}{2}\right) y+\frac{1}{2}y^T
	\mathbf{p}\left(\frac{\theta}{2}+\pi\right)^{H}\mathbf{f}\left(\frac{\theta}{2}+\pi\right)\mathbf{p}\left(\frac{\theta}{2}+\pi\right)y.\]\end{small}
	\begin{sloppypar}From the proof of the second item we already know that $y^T\hat{\mathbf{f}}y\ge 0,$ since both $\frac{1}{2}y^T\mathbf{p}\left(\frac{\theta}{2}\right)^{H}\mathbf{f}\left(\frac{\theta}{2}\right)\mathbf{p}\left(\frac{\theta}{2}\right) y$ and $\frac{1}{2}y^T
	\mathbf{p}\left(\frac{\theta}{2}+\pi\right)^{H}\mathbf{f}\left(\frac{\theta}{2}+\pi\right)\mathbf{p}\left(\frac{\theta}{2}+\pi\right)y$ are non negative. To prove that, for $\theta\neq 2\theta_0$, $y^T\hat{\mathbf{f}}y> 0,$ it is sufficient to show that if  $\frac{1}{2}y^T\mathbf{p}\left(\frac{\theta}{2}\right)^{H}\mathbf{f}\left(\frac{\theta}{2}\right)\mathbf{p}\left(\frac{\theta}{2}\right) y=0$, then $\frac{1}{2}y^T
	\mathbf{p}\left(\frac{\theta}{2}+\pi\right)^{H}\mathbf{f}\left(\frac{\theta}{2}+\pi\right)\mathbf{p}\left(\frac{\theta}{2}+\pi\right)y$ is different from 0 and, if $\frac{1}{2}y^T
	\mathbf{p}\left(\frac{\theta}{2}+\pi\right)^{H}\mathbf{f}\left(\frac{\theta}{2}+\pi\right)\mathbf{p}\left(\frac{\theta}{2}+\pi\right)y=0$, then $\frac{1}{2}y^T\mathbf{p}\left(\frac{\theta}{2}\right)^{H}\mathbf{f}\left(\frac{\theta}{2}\right)\mathbf{p}\left(\frac{\theta}{2}\right) y\neq 0$.
	In the following we show the first fact, the latter can be proved in the same way. 
	 We have that $\frac{1}{2}y^T\mathbf{p}\left(\frac{\theta}{2}\right)^{H}\mathbf{f}\left(\frac{\theta}{2}\right)\mathbf{p}\left(\frac{\theta}{2}\right) y=0$ if and only if $\mathbf{p}\left(\frac{\theta}{2}\right) y=0y$. Note that for $\theta\neq 2\theta_0$, $\mathbf{f}(\theta/2)>0$, since, for	$\theta\neq \theta_0$, $\mathbf{f}(\theta)>0$. Moreover,
	 if $\theta\neq 2 \theta_0$, then for periodicity, $\theta\neq 2 \theta_0-2\pi$, that is, $\theta/2+\pi\neq  \theta_0$. Consequently, $\mathbf{f}(\theta/2+\pi)>0$. If, by contradiction,  $\frac{1}{2}y^T
	\mathbf{p}\left(\frac{\theta}{2}+\pi\right)^{H}\mathbf{f}\left(\frac{\theta}{2}+\pi\right)\mathbf{p}\left(\frac{\theta}{2}+\pi\right)y= 0$, then  $\mathbf{p}(\theta/2+\pi)y=0y$, that cannot be satisfied, since, by hypothesis $\mathbf{p}(\theta)$ verifies the condition $(i)$.\end{sloppypar}
\item Since $\hat{\mathbf{f}}(\theta)$ is Hermitian for every $\theta$, we can assume $\hat{\mathbf{f}}(\theta)$ equals to
\begin{small}
\begin{equation*}\label{eq:dec_f_hat}
\begin{split}
&\left[
\begin{array}{@{\;}c@{\;}|@{\;}c@{\;}|@{\;}c@{\;}|@{\;}c@{\;}|@{\;}c@{\;}}
\!v_1(\theta)&\dots &v_{\bar{\jmath}}(\theta)&\dots& v_d(\theta)\!
\end{array}
\right]\!\!\begin{bmatrix}
\lambda_1(\hat{\mathbf{f}}(\theta))& & &  &\\
& \!\!\ddots\!\! & &  &\\
& & \lambda_{\bar{\jmath}}(\hat{\mathbf{f}}(\theta))& & \\
& & & \!\!\ddots\!\! & \\
& & & &\lambda_d(\hat{\mathbf{f}}(\theta))
\end{bmatrix}\!\!
\left[
\begin{array}{ccccccc}
{v_1}^H(\theta) \\
\hline
{v_2}^H(\theta)\\
\hline
\vdots\\
\hline
v_{\bar{\jmath}}^H(\theta)\\
\hline
\vdots\\
\hline
{v_d}^H(\theta)
\end{array}
\right].
\end{split}
\end{equation*}
\end{small}
From property 3, we can assume, without loss of generality, that 
\begin{equation}\label{eq:lim_eigenvec_hat_f}
\lim_{\theta\rightarrow \theta_0} v_{\bar{\jmath}}(2\theta)=q_{\bar{\jmath}}(\theta_0).
\end{equation}
Consequently, we obtain
\begin{small}
\begin{equation}
\begin{split}
&\lim_{\theta\rightarrow\theta_0}\!\frac{\lambda_{\bar{\jmath}}(\hat{\mathbf{f}}(2\theta))}{\lambda_{\bar{\jmath}}(\mathbf{f}(\theta))}\!=\!\lim_{\theta\rightarrow\theta_0}\!\frac{v_{\bar{\jmath}}(2\theta)^H \hat{\mathbf{f}}(2\theta)) v_{\bar{\jmath}}(2\theta)}{q_{\bar{\jmath}}(\theta) {\mathbf{f}}(\theta) q_{\bar{\jmath}}(\theta)}\!=\!\frac{1}{2} \lim_{\theta\rightarrow\theta_0}\!\!\left(\!\frac{v_{\bar{\jmath}}(2\theta)^H {\mathbf{p}}(\theta)^H\mathbf{f}(\theta)\mathbf{p}(\theta) v_{\bar{\jmath}}(2\theta)}{q_{\bar{\jmath}}(\theta) {\mathbf{f}}(\theta) q_{\bar{\jmath}}(\theta)}
+\right.\\
&\left.
\frac{v_{\bar{\jmath}}(2\theta)^H {\mathbf{p}}(\theta+\pi)^H\mathbf{f}(\theta+\pi)\mathbf{p}(\theta+\pi) v_{\bar{\jmath}}(2\theta)}{q_{\bar{\jmath}}(\theta) {\mathbf{f}}(\theta) q_{\bar{\jmath}}(\theta)}\right)= \frac{1}{2} \frac{q_{\bar{\jmath}}(\theta_0)^H {\mathbf{p}}(\theta_0)^H\mathbf{f}(\theta_0)\mathbf{p}(\theta_0) q_{\bar{\jmath}}(\theta_0)}{q_{\bar{\jmath}}(\theta_0) {\mathbf{f}}(\theta_0) q_{\bar{\jmath}}(\theta_0)}+ \\
&\frac{1}{2}
\lim_{\theta\rightarrow\theta_0} \frac{v_{\bar{\jmath}}(2\theta)^H {\mathbf{p}}(\theta+\pi)^H\mathbf{f}(\theta+\pi)\mathbf{p}(\theta+\pi) v_{\bar{\jmath}}(2\theta)}{q_{\bar{\jmath}}(\theta) {\mathbf{f}}(\theta) q_{\bar{\jmath}}(\theta)},
\end{split}
\end{equation}
\end{small}
where in the first term we used relation (\ref{eq:lim_eigenvec_hat_f}).  Exploiting properties (a) and (b) and Remark \ref{rmk:eig_p_pi}, we can write
\begin{equation}
\begin{split}
&\lim_{\theta\rightarrow\theta_0}\frac{\lambda_{\bar{\jmath}}(\hat{\mathbf{f}}(2\theta))}{\lambda_{\bar{\jmath}}(\mathbf{f}(\theta))}=\frac{1}{2}(\lambda^{(1)})^2+\frac{1}{2}{q_{\bar{\jmath}}(\theta_0)^H\mathbf{f}(\theta_0+\pi)q_{\bar{\jmath}}(\theta_0)}\left(\lim_{\theta\rightarrow\theta_0} \frac{\lambda_{\bar{\jmath}}^2(\mathbf{p}(\theta +\pi))}{\lambda_{\bar{\jmath}}(\mathbf{f}(\theta))}\right),
\end{split}
\end{equation}
where $\lambda_{\bar{\jmath}}^2(\mathbf{p}(\theta +\pi))$ is the eigenvalue of $\mathbf{p}(\theta +\pi)$ associated to the eigenvector $q_{\bar{\jmath}}(\theta_0)$. Since the quantity $\frac{1}{2}{q_{\bar{\jmath}}(\theta_0)^H\mathbf{f}(\theta_0+\pi)q_{\bar{\jmath}}(\theta_0)}$ is strictly positive and bounded by condition (\ref{eqn:condition_on_f}), property 5 is satisfied under the hypotheses that 
\begin{equation}
\lim_{\theta\rightarrow\theta_0} \frac{\lambda_{\bar{\jmath}}^2(\mathbf{p}(\theta +\pi))}{\lambda_{\bar{\jmath}}(\mathbf{f}(\theta))}<\infty.
\end{equation}

\end{enumerate}
\end{proof}
\end{lemma}  
{
Passing from TGM to V-cycle, the hypothesis 2 in Lemma \ref{lemm:condition_iii} has to be strengthen removing the power two, similarly to the case of scalar symbol \cite{AD}, obtaining the condition 
$$\lim_{\theta\rightarrow 	\theta_0} \frac{|\lambda_{\bar{\jmath}}(\mathbf{p}(\theta+\pi))|}{\lambda_{\bar{\jmath}}(\mathbf{f}(\theta))}=c.$$ 
In fact, it leads to convergence and optimality even when dealing with V-cycle with more than two grids.
}
\begin{lemma}\label{lem:convergence_vcycle}
Let $\mathbf{p}$ and $\mathbf{f}$ satisfy the hypothesis of Lemma \ref{lemma:properties_f_hat_f}. If 
\begin{equation}\label{eq:final_bound_vcycle}
\lim_{\theta\rightarrow 	\theta_0} \frac{|\lambda_{\bar{\jmath}}(\mathbf{p}(\theta+\pi))|}{\lambda_{\bar{\jmath}}(\mathbf{f}(\theta))}=c,
\end{equation}
then the two bounds in (\ref{eq:bound_1_2_vcycle}), needed for the  convergence and optimality of the V-cycle, are verified.
\begin{proof}
The condition  (\ref{eq:final_bound_vcycle}) implies the hypothesis of the item 5 of Lemma \ref{lemma:properties_f_hat_f}, then we have that the order of the zero at the coarse levels does not change, since it brings to 
\begin{equation}\label{eq:limit_f_fhat_in_proof}
\lim_{\theta\rightarrow\theta_0}\frac{\lambda_{\bar{\jmath}}(\hat{\mathbf{f}}(2\theta))}{\lambda_{\bar{\jmath}}(\mathbf{f}(\theta))}=c, \, c\neq 0\in \mathbb{R},
\end{equation}
where $ \hat{\mathbf{f}}(\theta)$ is given in (\ref{f2t}).  Hence, from direct computation using the same techniques as in the proof of Theorem \ref{th:tgmopt} and Lemma \ref{lemma:properties_f_hat_f}, we have that the  quantity 
$$\|\mathbf{f}(\theta) \mathbf{p}(\theta)\hat{\mathbf{f}}^{-1}(2\theta)\|$$
is bounded. Moreover, the limit bound in (\ref{eq:limit_f_fhat_in_proof}) implies that the second condition in (\ref{eq:bound_1_2_vcycle}) can be replaced by
$$\| \mathbf{p}(\theta+\pi){\mathbf{f}}^{-1}(2\theta)\|<\infty,$$
which is given by (\ref{eq:final_bound_vcycle}).
\end{proof}
\end{lemma}
\section{Geometric Projectors}
\label{sect:natural_projectors}
{In the following we will apply the theoretical considerations from the previous sections to problems arising from the discretization of partial differential equations (PDEs). When PDEs are discretized with high order of accuracy using the finite element method (FEM), block matrices arise naturally. We consider finite elements with nodal bases, using a \textcolor{red}{Cartesian} grid and Kronecker products of one-dimensional basis functions. For a problem posed in $m$ dimensions discretized using Kronecker product of basis functions of degree $r$ this automatically yields blocks of size $r^m \times r^m$. As prolongation operators we have different choices, here we consider two: the linear interpolation usually used in geometric multigrid methods for scalar problems \cite{Trot} and the prolongation obtained as the adjoint operator of the restriction operator when considering the finite element basis functions \cite{Brae}.}

\label{sect:fem_stiffness_matrix}
\subsection{$\mathbb{Q}_{r} $  Lagrangian FEM Stiffness Matrices} \label{sub:stiff_matrices}
{First we} consider the $\mathbb{Q}_{r} $ Lagrangian FEM approximation of the differential 1D problem{, that is given by:}
Find $u$ such that
\begin{small}
\begin{equation}\label{FEM_problem}
\begin{cases}
&-u''(x)=\psi(x) \quad {\rm on}\, \, (0,1), \\
& u(0)=u(1)=0 ,
\end{cases}
\end{equation}
\end{small}where $\psi(x)\in L^2\left(0,1\right)$. In this setting the weak formulation on the problem is written as follows:
Find $u\in H^{1}_0(0,1)$ such that
$
a(u,v)=\langle\psi,v\rangle, \quad \forall v\in H^1_0(0,1),$
where $a(u,v):=\int_{(0,1)}u'(x)v'(x)\, dx$ and  $\langle\psi,v\rangle := \int_{(0,1)}\psi(x)v(x)\, dx$.
For $r, n\ge 1$, we define the space 
\begin{small}
\begin{equation}
V_n^{(r)}:=\{\sigma\in \mathcal{C}\left([0,1]\right), \, {\rm such}\, {\rm that}\, \sigma{|_{\left[\frac{i}{n},\frac{i+1}{n}\right]}}\in \mathbb{P}_{r}, \, \forall i=0,\dots,n-1\},
\end{equation}
\end{small}where we denote by $\mathbb{P}_{r}$ the space of polynomials of degree less than or equal to $r$. {So} the space $V_n^{(r)}$ represents the space of continuous piecewise polynomial functions. Starting from $V_n^{(r)}$, we consider its subspace of functions that vanish on the boundary, defined by $
W_n^{(r)}:=\{\sigma\in V_n^{(r)}, \, {\rm such}\, {\rm that}\, \sigma(0)=\sigma(1)=0\}.$
Note that $W_n^{(r)}$ is a finite $n{r}-1$ dimension subspace of $H_0^1(0,1)$ and, following a Galerkin approach {\cite{Brae}}, we approximate the solution $u$ of the variational problem by solving the problem:
Find $u_{r,n}\in W_{n}^{(r)}$ such that
\begin{small}
\begin{equation}\label{eq:variational_problem}
a(u_{r,n},v)=\langle\psi,v\rangle, \quad \forall v\in W_n^{(r)}.
\end{equation}
\end{small}
We define the uniform knot sequence
\begin{small}
\begin{equation}\label{eq:knot_seq}
\xi_i=\frac{i}{n r}, \quad i=0,\dots, nr,
\end{equation}
\end{small}and the Lagrangian basis functions by $
\phi_j^{n,r}(\xi_i)=\delta_{i,j},\, i,j=0,\dots, nr,$
 with $\delta_{i,j}$ being the Kronecker delta.
{It is well known} that the latter definition is well-posed and that $\{\phi_1^{n,r},\dots,\phi_{nr-1}^{n,r}\}$ is a basis for $W_n^{(r)}$. Then $u_{r,n}$ can be written as {a} linear combination of such basis as
$
 u_{r,n}=\sum_{j=1}^{nr-1} u_j \phi_j^{n,r}{.}
$
 {Using this discretization, approximately solving} the problem (\ref{eq:variational_problem}) reduces to the {solution} of the linear system 
 $
 K_n^{(r)} \textbf{u}=\textbf{b},
$
 with 
 \begin{equation*}
 K_n^{(r)}=[a(\phi_j^{n,r},\phi_i^{n,r})]_{i,j=1}^{nr-1}, \quad \textbf{b}=\left[ \langle\psi,\phi_i^{n,r}\rangle\right]_{i=1}^{nr-1}, \quad \textbf{u}=\left[u_i\right]_{i=1}^{nr-1}.
 \end{equation*}
 The spectral properties of the \textcolor{red}{stiffness} matrix-sequence $\{ K_n^{(r)}\}_n$ were studied in \cite{qp}. In the following we report the spectral properties of matrix-valued function $\mathbf{f}$ associated with the normalized matrix-sequence $\{T_n(\mathbf{f})\}_n=\{n K_n^{(r)}\}_n,$  which are needed for our analysis \cite{qp,RST}.
 \begin{theorem}\label{thm:order_zero}
\begin{sloppypar}The $r \times r$ matrix-valued generating function of 
$\{T_n(\mathbf{f})\}_n=\{n K_n^{(r)}\}_n$ is 
\begin{align}
\mathbf{f}_{\mathbb{Q}_{r}}(\theta)=a_0+a_1{\rm e}^{\imath \theta}+a_1^{T}{\rm e}^{-\imath \theta}
\end{align}
and the following statements hold true:\end{sloppypar}
\begin{enumerate}
\item $\mathbf{f}_{\mathbb{Q}_{r}}(0) {\rm e}_{r}= 0$, 
${\rm e}_{r}$ vector of all ones, $r\ge 1$; 
\item there exist constants
$C_2\ge C_1>0$ (dependent on $f_{\mathbb{Q}_{r}}$) such that
\begin{equation}
C_1 (2-2\cos(\theta)) \le \lambda_1(\mathbf{f}_{\mathbb{Q}_{r}}(\theta))
\le C_2  (2-2\cos(\theta)); 
\end{equation}
 and 
\item there exist constants {$M_2\ge M_1>0$} (dependent on $\mathbf{f}_{\mathbb{Q}_{r}}$) such that
\begin{equation}
0 < { M_1} \le \lambda_j(\mathbf{f}_{\mathbb{Q}_{r}}(\theta)) \le {M_2},\ \ \ \ j=2,\ldots,r.
\end{equation}
\end{enumerate}
 \end{theorem}
 From the latter result, we have that $\mathbf{f}_{\mathbb{Q}_{r}}$ is a matrix-valued trigonometric polynomial which {fulfills} the hypotheses of Subsection \ref{sec:Theorem}.
  Indeed, for each $r\ge 1$ we have that 
  \begin{itemize}
  \item the minimum eigenvalue function $\lambda_{\min}(\mathbf{f}_{\mathbb{Q}_{r}})$ of $\mathbf{f}_{\mathbb{Q}_{r}}$ has a zero of order 2 
  in $\theta_0=0$;
  \item for $j=2,\dots, r$, it holds $\lambda_{j}(\mathbf{f}_{\mathbb{Q}_{r}}(\theta))>0$,
   for all $\theta\in [0,2\pi]$.
  \end{itemize} 
  Then, {the first item} of Theorem \ref{thm:order_zero} implies that 
   $\mathbf{f}_{\mathbb{Q}_{r}}$ can be decomposed as in equation (\ref{eq:dec_f}), with $q_{\bar{\jmath}}=q_1$ equal to ${\rm e}_{r}$, the column vector of all ones.  
   
   Consequently, in next section we test the applicability of the results given in \textcolor{red}{the} Section \ref{sec:Theorem}, and we confirm that two standard projectors are effective in terms of convergence and optimality when used for solving the linear system which has $T_n(\mathbf{f}_{\mathbb{Q}_{r}})$ as coefficient matrix.  In particular, we will deal with projector of the form 
   \begin{equation*}
P_{n,k}^{_{\mathbb{Q}_{r}}}=T_n(\mathbf{p}_{_{\mathbb{Q}_{r}  }})(K_{n,k}^{Even} \otimes I_d),
\end{equation*}
where its properties and efficiency will depend on those of the associated trigonometric polynomial $\mathbf{p}_{_{\mathbb{Q}_{r}  }}$. 

\subsection{The scalar linear interpolation projector}\label{sec:linear_interpolation}
{
The discretization using $\mathbb{Q}_r$ Lagrangian finite elements provides an approximate solution at all nodes used as interpolation points within an element and on its boundaries.
While the use of linear interpolation is common in finite difference discretizations of partial differential equations \cite{Trot}, it can be used in this setting, as well.
For $n$ elements and using polynomial degree $r$ it can be written as
}%
\begin{equation}\label{eqn:linear_interp_p}
P_{n,k}^{{r}}=T_{r n} (2+2\cos(\theta)) K^{Even}_{rn, \frac{r(n-1)}{2}}.
\end{equation}
{We will now show that this scalar interpolation nevertheless satisfies the} hypotheses of Lemma \ref{lemm:condition_on_s}. {For that purpose} we have to show that it fits into the block setting of the present paper. Hence, first we have to rewrite $P_{n,k}^{r}$ in block form as
\begin{equation}\label{eqn:linear_interp_p_blocked}
P_{n,k}^{{r}}=T_{n}(\mathbf{p}_{_{L_{r}}})(K^{Even}_{n,k} \otimes I_{r}).
\end{equation}
which means that we want to find a matrix-valued trigonometric polynomial $\mathbf{p}_{_{L_{r}}}$ such that the latter equation is true, with $P_{n,k}^{{r}}$ defined as in (\ref{eqn:linear_interp_p}).

Recalling the action of the cutting matrix $(K^{Even}_{n,k} \otimes I_{r})$, seen in Subsection \ref{sect:project_circ}, we {observe} that $P_{n,k}^{{r} }$ can be rewritten in the desired block form with associated matrix-valued trigonometric polynomial $\mathbf{p}_{_{L_{r}}}$ of the form
\begin{equation}\label{eq:pol_bolten}
\mathbf{p}_{_{L_{r}}}(\theta) = \hat{a}_0+\hat{a}_{-1} {\rm e}^{-\imath\theta}+\hat{a}_{1} {\rm e}^{\imath\theta},
\end{equation}
where the expression of the Fourier coefficients  $\hat{a}_0,\hat{a}_{-1}, \hat{a}_{1} $ depends on whether the degree is even or odd.
Indeed, we have
{
\begin{enumerate}
 \item In the case of even degree $r$, we define 
  \begin{equation*}
   A_{1}= T_{r} (2+2\cos(\theta)) K^{Even}_{r, \frac{r}{2}}=\begin{bmatrix}
    1 & & & & \\
    2 & & & & \\
    1 & 1 & & & \\
    & 2 & & & \\
    & 1 & & & \\
    & & \ddots & & \\
    & & & 1 & \\
    & & & 2 & \\
    & & & 1 & 1 \\
    & & & & 2
   \end{bmatrix} \in \mathbb{R}^{r\times \frac{r}{2}}.
  \end{equation*}
  Then the identity \eqref{eqn:linear_interp_p_blocked} holds with
  \begin{align*}
 \hat{a}_{-1} & = \begin{bmatrix}
    & A_1 & | & \textbf{0}_{{r,\frac{r}{2}}} & 
   \end{bmatrix}, \quad   \hat{a}_0 = \left[\begin{array}{ccc|c|c}
    & \textbf{0}_{1,\frac{r}{2}-1} & & 1 &  \\
    \cline{1-4}
    & & & & A_1\\
    & \textbf{0}_{{r}-1,\frac{r}{2}-1} & & \textbf{0}^T_{1,\frac{r}{2}-1} &
   \end{array}\right], \\
   \hat{a}_{1} & = \left[\begin{array}{ccc|c}
    & \textbf{0}_{1,r-1} & & 1\\
    \cline{1-4} & & & \\
    & \textbf{0}_{r-1,r-1} & & \textbf{0}^T_{r-1}
   \end{array}\right].
  \end{align*}
  Note that
  \begin{equation}\label{eq:structure_A1}
   \sum_{j=1}^{\frac{r}{2}} [A_1]_{1,j}=1 \quad \text{ and } \quad \sum_{j=1}^{\frac{r}{2}} [A_1]_{i,j}=2,\, \text{ for $i=2,\dots,r$.}
  \end{equation}
 \item In the case of odd degree $r$ we define
 \begin{small}
  \begin{align*}
   A_2 = T_r(2+2\cos(\theta)) K^{Odd}_{r,\frac{r+1}{2}} & = \begin{bmatrix}
    2 & & & & \\
    1 & 1 & & & \\
    & 2 & & & \\
    & 1 & & & \\
    & & \ddots & & \\
    & & & 1 & \\
    & & & 2 & \\
    & & & 1 & 1 \\
    & & & & 2
   \end{bmatrix} \in \mathbb{R}^{r\times \frac{r+1}{2}},\\
   A_3 = T_{r} (2+2\cos(\theta)) K^{Even}_{r, \frac{r+1}{2}} & = \begin{bmatrix}
    1 & & & & \\
    2 & & & & \\
    1 & 1 & & & \\
    & 2 & & & \\
    & 1 & & & \\
    & & \ddots & & \\
    & & & 1 & \\
    & & & 2 & \\
    & & & 1 & 1
   \end{bmatrix} \in \mathbb{R}^{r\times \frac{r+1}{2}}.
  \end{align*}
 \end{small}
  Then \eqref{eqn:linear_interp_p_blocked} holds with
 \begin{small}
  \begin{align*}
   \hat{a}_{-1} & = \begin{bmatrix}
    A_3 & | & \textbf{0}_{r,\frac{r-1}{2}} & 
   \end{bmatrix}, \\
   \hat{a}_0 & = \begin{bmatrix}
    & \textbf{0}_{r,\frac{r-1}{2}} & | &	A_2 &
   \end{bmatrix}, \\
   \hat{a}_{1} & = \left[\begin{array}{ccc|c}
    & \textbf{0}_{1,r-1} & & 1\\
    \cline{1-4}
    & & & \\
    & \textbf{0}_{r-1,r-1} & & \textbf{0}^T_{r-1}
   \end{array}\right].
  \end{align*}
 \end{small}
  Further we have
  \begin{small}
  \begin{equation}\label{eq:structure_A2_A3}
   \sum_{j=1}^{\frac{r+1}{2}} [A_2]_{1,j}=2, \, \sum_{j=1}^{\frac{r+1}{2}} [A_3]_{1,j}=1 \, \text{ and } \, \sum_{j=1}^{\frac{r+1}{2}} [A_2]_{i,j}=\sum_{j=1}^{\frac{r+1}{2}} [A_3]_{i,j}=2 \text{ for $i=2,\dots,r$.}
  \end{equation}
  \end{small}
\end{enumerate}
}%
In the following we show that $\mathbf{p}_{_{L_{r}}}$ satisfies the hypotheses of the Lemma \ref{lemm:condition_on_s}.
\begin{lemma}\label{lemma:verify_p_bolten}
Let $\mathbf{p}_{_{L_{r}}}$ be the $r\times r$ trigonometric polynomial defined in (\ref{eq:pol_bolten}), and ${\rm e}_{r}=[1,\dots,1]^T\in \mathbb{R}^r$. Then,
 \begin{enumerate}
 \item [1)]  $\mathbf{p}_{_{L_{r}}}(0)\,{\rm e}_{r}=4\, {\rm e}_{r}$.
 \item  [2)]  $\mathbf{p}_{_{L_{r}}}(\pi)\,{\rm e}_{r}=0\, {\rm e}_{r} $.
   \item  [3)]  $\mathbf{p}_{_{L_{r}}}(0)^H\,{\rm e}_{r}=4\, {\rm e}_{r}$. 
 \end{enumerate}
\begin{proof}
The first two items are equivalent to require that the sum of the elements in each row of the matrices $\mathbf{p}_{_{L_{r}}}(0)$  and $\mathbf{p}_{_{L_{r}}}(\pi)$ is $4$ and $0$, respectively.
 Hence, to prove $1)$, it is sufficient to show that
 \begin{equation*}
 \sum_{j=1}^{r} [\mathbf{p}_{_{L_{r}}}(0)]_{i,j}=\sum_{j=1}^{r}[\hat{a}_0+\hat{a}_1+\hat{a}_{-1}]_{i,j}=4, \quad i=1,\dots, r.
 \end{equation*}
Then, we can exploit the structure of the Fourier coefficients $\hat{a}_{-1}, \hat{a}_{1}, \hat{a}_{0}   $ for even and odd degree. In particular, looking at the structure of the matrices $A_1$, $A_2$, $A_3$ and at relations {\eqref{eq:structure_A1} and \eqref{eq:structure_A2_A3}}, we have that for even degree $r$ 
\begin{small} 
 \begin{equation*}
 \begin{split}
& \sum_{j=1}^{r} [\mathbf{p}_{_{L_{r}}}(0)]_{i,j}= \begin{cases}
 1+\left( 2\sum_{j=1}^{\frac{r}{2}}[A_1]_{1,j}\right)+1=4,  &{\rm for} \, i=1\\
\left( 2\sum_{j=1}^{\frac{r}{2}}[A_1]_{i,j}\right)=4, 	  &{\rm for} \, i=2,\dots,r\\
\end{cases},
 \end{split}
 \end{equation*}
 \end{small}
 and for odd degree $r$
 \begin{small}
 \begin{equation*}
 \begin{split}
& \sum_{j=1}^{r} [\mathbf{p}_{_{L_{r}}}(0)]_{i,j}\!=	\begin{cases}
\left( \sum_{j=1}^{\frac{r+1}{2}}[A_3]_{1,j}\!+\![A_2]_{1,j}\right)\!+\!1\!=\!4, &	{\rm for} \, i\!=\!1\\
\left( \sum_{j=1}^{\frac{r}{2}}[A_3]_{i,j}\!+\![A_2]_{i,j}\right)\!=\!4, & {\rm for} \, i\!=\!2,\dots,r
\end{cases}\!.
 \end{split}
 \end{equation*}
 \end{small}
The proof of $2)$ can be repeated following the idea in $1)$ and noting that
\[\mathbf{p}_{_{L_{r}}}(\pi)=\hat{a}_0-\hat{a}_1-\hat{a}_{-1}.\] 
Analogously, the third item can be proven following the same idea of $1)$, showing that the sum of the elements in each column  of the matrices $\mathbf{p}_{_{L_{r}}}(0)$ is $4$. Since it is a straightforward computation, we omit the details.
\end{proof}
 \end{lemma}
The latter result, together with Lemma \ref{lemm:condition_on_s} permits to conclude that $\mathbf{p}_{_{L_{r}}}$ satisfies condition $(ii)$, once that we prove that it satisfies the condition $(i)$, 
so that the matrix-valued function $\textbf{s}$ is well-defined.
By direct computation,
we find that for both even and odd $r$ we have
 $$\mathbf{p}_{_{L_{r}}}(\theta)^{H}\mathbf{p}_{_{L_{r}}}(\theta)+\mathbf{p}_{_{L_{r}}}(\theta+\pi)^{H}\mathbf{p}_{_{L_{r}}}(\theta+\pi)=\begin{bmatrix}
12&2& 0& \dots &2{\rm e}^{2\imath\theta}\\
2 & 12& 2&\dots &0\\ 
& &\ddots  & \\
0 & & & 12&2\\
2{\rm e}^{-2\imath\theta} & 0 &\dots&2 &12
\end{bmatrix},$$
which is clearly a definite positive matrix for all $\theta\in[0,2\pi)$, so $\mathbf{p}_{_{L_{r}}}$ satisfies condition $(i)$.
 Then, the function $ \textbf{s}(\theta)$ defined in (\ref{eqn:def_s}) is well-defined.

{
The validation of condition $(iii)$ will be investigated in subsection \ref{sec:V-cycle_Qp} since we will treat it for both the projector associated with $\mathbf{p}_{_{L_{r}}}$ and the classical geometric projector described in the next subsection.
}

\subsection{{The projector using the finite element basis functions}} \label{sec:geometric_projector}
In the present subsection we deal with the projector which is constructed following the {classical approach used for finite elements \cite{Brae}}. It has been already treated for the $\mathbb{Q}_{r} $ Lagrangian FEM stiffness matrices in the algebraic multigrid setting. Indeed, in \cite{FRTBS} the authors proved the optimality of the TGM methods for $r=2,3$.
Our goal is to generalize the study of the multigrid method procedure, proving that the restriction matrix can be written for any degree $r$ in the form 
\begin{equation}\label{eqn:form_project}
P_{n,k}^{{r}}=T_n(\mathbf{p}_{_{G_{r}  }})(K_{n,k}^{Even} \otimes I_d),
\end{equation}
with $\mathbf{p}_{_{G_{r}  }}$ being a matrix-valued trigonometric polynomial that satisfies the hypotheses of Lemma \ref{lemm:condition_on_s}, and then those of Theorem \ref{th:tgmopt}.

Let us start fixing the $r$ and taking as first $n=2$. Then we have from (\ref{eq:knot_seq}) the knot sequence $
\xi_i^{2r}=\frac{i}{2 r}, \, i =0,\dots,2r.$ If we take $n=4$ the uniform knot sequence is 
$
\xi_i^{4r}=\frac{i}{4 r}, \, i =0,\dots,4r,
$
and \textcolor{red}{thus} it can be obtained from $\{\xi_i^{2r}\}_i$ adding the midpoint of each sub interval defined by the points in  $\{\xi_i^{2r}\}_i$.

Taking the Lagrangian basis \textcolor{red}{functions} for the spaces $W_2^{r}$ and $W_4^{r}$, defined in Subsection \ref{sub:stiff_matrices},  the latter observation implies that 
\begin{equation}\label{eq:inclusion_spaces}
W_2^{r}\subseteq W_4^{r}.
\end{equation}
\begin{sloppypar}The geometric multigrid strategy suggests to construct a prolongation operator 
$\mathcal{P}~:~W_2^{r}\rightarrow W_4^{r}$,
imposing
$
\mathcal{P}v^{2,r}=v^{2,r},$ $\forall v^{2,r}\in  W_2^{r}.
$
A basis function $\phi_j^{2,r}$ can be then written as a linear combination of the functions $\phi_j^{4,r}$, that is $
\phi_j^{2,r}(x)=\sum_{i=1}^{4r-1} c_i \phi_i^{4,r}(x).$
From the properties of the basis functions, we have that $\phi_j^{4,r}(\xi_i)=r_{i,j}$.
Consequently, for $i=1,\dots, 4r-1$, the coefficients $c_i$ are given by the evaluations of $\phi_j^{2,r}(\xi_i) $, and this implies that
$
\phi_j^{2,r}(x)=\sum_{i=1}^{4r-1}\phi_j^{2,r}(\xi_i)  \phi_i^{4,r}(x).
$
Therefore, the $j-th$ column $({P})_j$ of the matrix $P$ representing the prolongation operator $\mathcal{P}$ is given by 
\begin{equation}\label{eq:jcolumn_P}
({P})_j=[\phi_j^{2,r}(\xi_i)]_{i=1}^{4 r-1}.
\end{equation}\end{sloppypar}
Since we are in the setting of multigrid methods for $r\times r$ block-Toeplitz matrices,  from Subsection \ref{sect:project_circ}  we have that  $n$ is taken of the form $2^t-1$, $k=\frac{n-1}{2}$ and we look for a prolongation matrix of the form (\ref{eqn:form_project}).
Taking inspiration from \cite{FRTBS}, we define the matrix-valued trigonometric polynomial $\mathbf{p}_{_{G_{r}  }}$ by
\begin{equation}\label{eq:pol_cristina}
\mathbf{p}_{_{G_{r}  }}(\theta)=\hat{b}_0+ \hat{b}_1 {\rm e}^{\iota \theta}+ \hat{b}_{-1} {\rm e}^{-\iota \theta}+ \hat{b}_2 {\rm e}^{2\iota \theta}+ \hat{b}_{-2} {\rm e}^{-2\iota \theta}
\end{equation}
with
{$
\hat{b}_{-2}=0_{r\times r},\,
\hat{b}_{-1}=[P_{i,j}]_{\substack{i=1,\dots,r \\ j=1,\dots, r}},\,
\hat{b}_0=[P_{i,j}]_{\substack{i=r+1,\dots,2r \\ j=1,\dots, r}},\,
\hat{b}_1=[P_{i,j}]_{\substack{i=2r+1,\dots,3r \\ j=1,\dots, r}}, \,$ and
$\hat{b}_{2}=[P_{i,j}]_{\substack{i=3r+1,\dots,4r \\ j=1,\dots, r}}.$
}%
Then, from the expression of the columns of the matrix $P$ in (\ref{eq:jcolumn_P}), we have
{
\begin{align*}
\hat{b}_{-1}=\begin{bmatrix}
\phi_1^{2,r}(\xi_{1}^{4,r})&  \dots& \phi_{{\rm r}}^{2,r}(\xi_{1}^{4,r})\\
\vdots&  & \vdots\\
\phi_1^{2,r}(\xi_{r}^{4,r})&  \dots& \phi_{{\rm r}}^{2,r}(\xi_{r}^{4,r}) \\
\end{bmatrix}, \quad &
\hat{b}_0=\begin{bmatrix}
\phi_1^{2,r}(\xi_{r+1}^{4,r})&  \dots& \phi_{{\rm r}}^{2,r}(\xi_{r+1}^{4,r})\\
\vdots&  & \vdots\\
\phi_1^{2,r}(\xi_{2r}^{4,r})&  \dots& \phi_{{\rm r}}^{2,r}(\xi_{2r}^{4,r}) \\
\end{bmatrix}, \\
\hat{b}_{1}=\begin{bmatrix}
0&  \dots & 0 & \phi_{{\rm r}}^{2,r}(\xi_{2r+1}^{4,r})\\
\vdots &  & \vdots&\vdots\\
0&  \dots & 0 & \phi_{{\rm r}}^{2,r}(\xi_{3r}^{4,r})\\
\end{bmatrix}, \quad &
\hat{b}_{2}=\begin{bmatrix}
0&  \dots & 0 & \phi_{{\rm r}}^{2,r}(\xi_{3r+1}^{4,r})\\
\vdots &  & \vdots&\vdots\\
0&  \dots & 0 & \phi_{{\rm r}}^{2,r}(\xi_{4r}^{4,r})\\
\end{bmatrix}.
\end{align*}
}%
\begin{sloppy}where, for the expressions of  $\hat{b}_1$ and $\hat{b}_2$, we are using the fact that the sets
 ${\rm supp}(\phi_1^{2,r}),\dots, {\rm supp}(\phi_{r-1}^{2,r})$ are included in $[0,\xi_{2r}^{4,r}]$, see \cite[pag 1108]{qp}. Moreover, from \cite[Equation (3.5)]{qp} we can also see that, for $i=2r+1,\dots,4r-1$, 
 \begin{equation*}
 \phi_{r}^{2,r}(\xi_i^{4,r})=\phi_0^{2,r}\left(\xi_i^{4,r}-\frac{1}{2}\right)=\phi_0^{2,r}(\xi_{i-2r}^{4,r}),
 \end{equation*}\end{sloppy}
 hence we have
 \begin{equation}
 \begin{split}
\hat{b}_{1}=\begin{bmatrix}
0&  \dots & 0 & \phi_{0}^{2,r}(\xi_{1}^{4,r})\\
\vdots &  & \vdots&\vdots\\
0&  \dots & 0 & \phi_{0}^{2,r}(\xi_{r}^{4,r})\\
\end{bmatrix} \quad {\rm and} \quad \hat{b}_{2}=\begin{bmatrix}
0&  \dots & 0 & \phi_{0}^{2,r}(\xi_{r+1}^{4,r})\\
\vdots &  & \vdots&\vdots\\
0&  \dots & 0 & \phi_{0}^{2,r}(\xi_{2r}^{4,r})\\
\end{bmatrix}. \\
\end{split}
 \end{equation}
We want to prove that such projector $P^{{r}  }_{n,k}$ satisfies the hypothesis  of  Lemma \ref{lemm:condition_on_s}. 
Since, from Theorem \ref{thm:order_zero}, we know that ${\rm e}_{r}=[1,\dots,1]^T\in\mathbb{R}^r$ is the eigenvector of $\mathbf{f}(0)$ associated with the ill-conditioned subspace, the next lemma gives us the proof that $P^{{r}  }_{n,k}$ satisfies the hypothesis of Theorem \ref{th:tgmopt}.
 \begin{lemma}\label{lemma:verify_p_cristina}
Let $\mathbf{p}_{_{G_{r}  }}$ be the $r\times r$ trigonometric polynomial defined in (\ref{eq:pol_cristina}), and ${\rm e}_{r}=[1,\dots,1]^T\in \mathbb{R}^r$. Then 
 \begin{enumerate}
 \item [1)] $\mathbf{p}_{_{G_{r}  }}(0)\,{\rm e}_{r}=2\, {\rm e}_{r}$,
 \item [2)] $\mathbf{p}_{_{G_{r}  }}(\pi)\,{\rm e}_{r}=0\, {\rm e}_{r} $,
  \item [3)] {$\mathbf{p}_{_{G_{r}  }}(0) $ is non-singular.}
 \end{enumerate}
 \begin{proof}
 Note that the thesis is equivalent to require that the sum of the elements in each row of the matrices $\mathbf{p}_{_{G_{r}  }}(0)$  and $\mathbf{p}_{_{G_{r}  }}(\pi)$ is $2$ and $0$, respectively.
 Then, to prove item $1)$, we prove that for every $i=1,\dots, r$ $
 \sum_{j=1}^{r} [\mathbf{p}_{_{G_{r}  }}(0)]_{i,j}=2.$
The expression of $\mathbf{p}_{_{G_{r}  }}(\theta)$ in (\ref{eq:pol_cristina}) yields
$\mathbf{p}_{_{G_{r}  }}(0)=\hat{b}_0+\hat{b}_1+\hat{b}_{-1}+\hat{b}_{2},$ then we have for $i=1,\dots, r$
\begin{equation}
\begin{split}
&\sum_{j=1}^{r} [\mathbf{p}_{_{G_{r}  }}(0)]_{i,j}=\sum_{j=1}^{r} [\hat{b}_0+\hat{b}_1+\hat{b}_{-1}+\hat{b}_{2}]_{i,j}=\\
&\left( \sum_{j=1}^{r} \phi_j^{2,r}(\xi_{i+r}^{4,r})\right)+ \phi_0^{2,r}(\xi_{i}^{4,r})+\left( \sum_{j=1}^{r} \phi_j^{2,r}(\xi_{i}^{4,r})\right)+  \phi_{0}^{2,r}(\xi_{i+r}^{4,r})
\end{split}.
\end{equation}
\begin{sloppy}From the fact that the basis functions form a partition of the unity and since
${\rm supp}(\phi_j^{2,r})\cap [0,\xi^{4,r}_{2r}]=\emptyset$, $j=r+1,\dots,2r,$ we have, for $i=1,\dots, r$, that 
\begin{equation*}
\sum_{j=1}^{r} [\mathbf{p}_{_{G_{r}  }}(0)]_{i,j}= \left( \sum_{j=0}^{2r} \phi_j^{2,r}(\xi_{i}^{4,r})\right)+  \left( \sum_{j=0}^{2r} \phi_j^{2,r}(\xi_{i+r}^{4,r})\right)=2.
\end{equation*}\end{sloppy}
In order to prove item $2)$ we write analogously $\mathbf{p}_{_{G_{r}  }}(\pi)=\hat{b}_0-\hat{b}_1-\hat{b}_{-1}+\hat{b}_{2},$ and, for $i=1,\dots,r$,
\begin{equation*}
\begin{split}
&\sum_{j=1}^{r} [\mathbf{p}_{_{G_{r}  }}(\pi)]_{i,j}=\sum_{j=1}^{r} [\hat{b}_0\textcolor{black}{+}\hat{b}_1-\hat{b}_{-1}+\hat{b}_{2}]_{i,j}=\\
& \left( \sum_{j=0}^{2r} \phi_j^{2,r}(\xi_{i}^{4,r})\right)- \left( \sum_{j=0}^{2r} \phi_j^{2,r}(\xi_{i+r}^{4,r})\right)=0.
\end{split}
\end{equation*}
The item $3)$ is equivalent to  ${\rm det}\,\left(\mathbf{p}_{_{G_{r}  }}(0)\right)\neq0$. By direct computation, we have that $${\rm det}\,\left(\mathbf{p}_{_{G_{r}  }}(\theta)\right)=\frac{{\rm e}^{-r\iota \theta}({\rm e}^{\iota \theta}+1)^{r+1}}{2^{\frac{r(r+1)}{2}}}.$$
Finally, since ${\rm e}^{-r\iota \theta}\neq 0$, $\forall \theta\in [0,2\pi]$ and $({\rm e}^{\iota \theta}+1)=0$ only if $\theta=\pi+2\ell\pi$, we have that $${\rm det}\,\left(\mathbf{p}_{_{G_{r}  }}(\theta)\right)\neq0, \quad {\rm for}\quad \theta \in [0,2\pi]\setminus \{\pi\},$$
hence $\mathbf{p}_{_{G_{r}  }}(0)$ is non-singular. 
 \end{proof}
 \end{lemma}
Once we verify that $\mathbf{p}_{_{G_{r}}}$ satisfies condition $(i)$, we can use Lemmas \ref{lemma:verify_p_cristina} and \ref{lemm:condition_on_s_p_invertible} to conclude that the matrix-valued function $\textbf{s}$ is well-defined and $\mathbf{p}_{_{G_{r}}}$ satisfies condition $(ii)$.
To prove that 
$\mathbf{p}_{_{G_{r}}}(\theta)^{H}\mathbf{p}_{_{G_{r}}}(\theta)+\mathbf{p}_{_{G_{r}}}(\theta+\pi)^{H}\mathbf{p}_{_{G_{r}}}(\theta+\pi)>0$ it is sufficient to show that both ${\rm det}\,(\mathbf{p}_{_{G_{r}  }}(\theta)^H\mathbf{p}_{_{G_{r}  }}(\theta))$ and ${\rm det}\,(\mathbf{p}_{_{G_{r}  }}(\theta+\pi)^H\mathbf{p}_{_{G_{r}  }}(\theta+\pi))$ are non-negative definite \textcolor{red}{matrix-valued functions} which are singular respectively in $\theta_1$ and $\theta_2$, with $\theta_1\neq\theta_2.$ 
Indeed, we have that
 $${\rm det}\,\left(\mathbf{p}_{_{G_{r}  }}(\theta)^H\mathbf{p}_{_{G_{r}  }}(\theta)\right)={\rm det}\,(\mathbf{p}_{_{G_{r}  }}(\theta))^2=\frac{{\rm e}^{-2r\iota \theta}({\rm e}^{\iota \theta}+1)^{2(r+1)}}{2^{r(r+1)}},$$
 which is zero for $\theta_1=\pi.$
 Analogously, it holds
 \begin{equation*}
 \begin{split}
& {\rm det}\,\left(\mathbf{p}_{_{G_{r}  }}(\theta+\pi)^H\mathbf{p}_{_{G_{r}  }}(\theta+\pi)\right)=\\
 &{\rm det}\,\left(\mathbf{p}_{_{G_{r}  }}(\theta+\pi)\right)^2=\left(\frac{{\rm e}^{-r\iota (\theta+\pi)}({\rm e}^{\iota (\theta+\pi)}+1)^{r+1}}{2^{\frac{r(r+1)}{2}}}\right)^2=\frac{{\rm e}^{-2r\iota \theta}({\rm e}^{\iota \theta}-1)^{2(r+1)}}{2^{r(r+1)}},
 \end{split}
 \end{equation*}
  which is zero for $\theta_2=0.$

 \subsection{Optimal convergence of the V-cycle using the projector $P_{n,k}^{_{\mathbb{Q}_{r}}}$}
\label{sec:V-cycle_Qp}
 
{For both projectors described in Subsections \ref{sec:linear_interpolation} and \ref{sec:geometric_projector} we have to verify the limit condition {(iii)}, in order to theoretically ensure the TGM optimality.  }

{
For this purpose it is sufficient to show either that  the function  $1-\lambda_{\bar{\jmath}}(\textbf{s}(\theta))$ has a zero at least of the same order of $\lambda_{\bar{\jmath}}(\mathbf{f}_{\mathbb{Q}_{r}}(\theta))$, or, using the result in Lemma \ref{lemm:condition_iii}, that this property is satisfied by the eigenvalue function $\lambda_{\bar{\jmath}}(\mathbf{p}_{_{G_{r}}}(\theta+\pi))^2$.
}

\begin{enumerate}
\item {For the linear interpolation operator:}
{
\begin{itemize}
\item  for even degree, we have that $\textbf{s}(\theta)$ is a projector since it can be easily verified that $\textbf{s}^2(\theta)-\textbf{s}(\theta)=\textbf{0}_{r\times r}$. Hence, from condition {(ii)}, we have $\lambda_{\bar{\jmath}}(\textbf{s}(0))= 1$, and, from the continuity of the eigenvalue functions (Lemma \ref{lemm:continuity_eig}), we have that $\lambda_{\bar{\jmath}}(\textbf{s}(\theta))\equiv 1$. Hence, it is straightforward to see that the  condition {(iii)} is verified;
\item for odd degree, it can be numerically proved that $\lambda_{\bar{\jmath}}(\mathbf{p}_{_{L_{r}}}(\theta+\pi))^2$ has a zero of order $4$ in $0$, then  condition {(iii)} is verified.
For this purpose, we can numerically study the behavior of the function $\det(\mathbf{p}_{_{L_{r}}}(\theta+\pi)\mathbf{p}_{_{L_{r}}}(\theta+\pi)^H)$. Indeed, since for $l=2,\dots,r, \,  \theta \in [0,2\pi],$
\begin{equation}\label{eq:lambda_1vsall}
\lambda_{\bar{\jmath}}\left(\mathbf{p}_{_{L_{r}}}(\theta+\pi)\mathbf{p}_{_{L_{r}}}(\theta+\pi)^H\right)< \lambda_l \left(\mathbf{p}_{_{L_{r}}}(\theta+\pi)\mathbf{p}_{_{L_{r}}}(\theta+\pi)^H\right), 
\end{equation}
 the behaviour of $\lambda_{\bar{\jmath}}(\mathbf{p}_{_{L_{r}}}(\theta+\pi)\mathbf{p}_{_{L_{r}}}(\theta+\pi)^H)$ in $0$ is equivalent to that of $
\det\left(\mathbf{p}_{_{L_{ {r}}}}(\theta+\pi)\mathbf{p}_{_{L_{ {r}}}}(\theta+\pi)^H\right)= \prod_{i=1}^r \lambda_l \left(\mathbf{p}_{_{L_{ {r}}}}(\theta+\pi)\mathbf{p}_{_{L_{ {r}}}}(\theta+\pi)^H\right)$ 
at the same point, which as a product of nonnegative functions is still a nonnegative function. We numerically checked that
\begin{align*}
\det\left(\mathbf{p}_{_{L_{ {r}}}}(\theta+\pi)\mathbf{p}_{_{L_{ {r}}}}(\theta+\pi)^H\right)= {\rm e}^{-2\iota\theta}({\rm e}^{\iota\theta} - 1)^4,
\end{align*}
which has a zero of order $4$ in $0$.
\end{itemize}
}
\item { For the geometric projector operator we consider the even and odd degree $r$ simultaneously and we follow the latter strategy of studying the behavior of $\det\left(\mathbf{p}_{_{G_{ {r}}}}(\theta+\pi)\mathbf{p}_{_{G_{ {r}}}}(\theta+\pi)^H\right)$ in $0$. From the proof of item 3) of Lemma \ref{lemma:verify_p_cristina} we have that
\begin{align*}
\det\left(\mathbf{p}_{_{G_{ {r}}}}(\theta+\pi)\mathbf{p}_{_{G_{ {r}}}}(\theta+\pi)^H\right)=\frac{{\rm e}^{-2r\iota \theta}({\rm e}^{\iota \theta}-1)^{2(r+1)}}{2^{r(r+1)}},
\end{align*}  
which clearly has a zero of order $2(r+1)$ in $0$.}
\end{enumerate}
\textcolor{red}{
\subsection{Example $r=2$}
We conclude the section with the application of the results in Section \ref{sec:tgmoptimality} and Subsection \ref{sec:condition_Vcycle} on the scalar linear interpolation projector for the specific case $r=2$ for the problem
\begin{equation*}
 \left\{\begin{array}{ll}
	-(a(x)u(x)')' = \psi(x),\, {\rm on} \, (0,1)\\	
	u(0)=u(1)=0,
	\end{array}\right.
\end{equation*}
which is the variable coefficient version of the problem in (\ref{FEM_problem}). 
In this setting the grid transfer operator $P_{n,k}^{{2}}$ is a $2n\times n-1$ matrix given by
}
\textcolor{red}{
\begin{equation}
P_{n,k}^{{2}}=T_{n}(\mathbf{p}_{_{L_{2}}})(K^{Even}_{n,k} \otimes I_{2}).
\end{equation}
with matrix-valued trigonometric polynomial $\mathbf{p}_{_{L_{2}}}(\theta)=\begin{bmatrix}
1 +{\rm e}^{-\iota \theta}& {\rm e}^{\iota \theta}+1\\
2{\rm e}^{-\iota \theta} & 2
\end{bmatrix}.$ 
}
\textcolor{red}{
The cut stiffness matrix 	$[T_n(\mathbf{\textbf{f}}_{\mathbb{Q}_2})]_-,$ has the associated generating function  
\begin{equation*}
\begin{split}
\mathbf{f}_{\mathbb{Q}_2}(\theta)=&\frac{1}{3}\left(\begin{bmatrix}
16 & -8 \\
-8 &14
\end{bmatrix}+ \begin{bmatrix}
0 & -8 \\
0 &1
\end{bmatrix}{\rm e}^{\iota\theta}+\begin{bmatrix}
0 & 0\\
-8 &1
\end{bmatrix} {\rm e}^{-{\iota}				\theta}\right)=\\&\frac{1}{3}\begin{bmatrix}
16 & -8(1+{\rm e}^{\iota \theta})\\
-8(1+{\rm e}^{-\iota \theta})& 14+ {\rm e}^{\iota \theta}+{\rm e}^{-\iota \theta}\end{bmatrix},
\end{split}
\end{equation*}
with the following properties:
		\begin{itemize}
		\item $\lambda_{1}(\mathbf{\textbf{f}}_{\mathbb{Q}_2}(\theta))$ has a zero of order 2 in $\theta_0=0$.
		\item $\mathbf{\textbf{f}}_{\mathbb{Q}_{{2}}}(0) q_{1}(0)= 0$, $q_{1}(0)=[1,1]^T$.
		\end{itemize}
}
\textcolor{red}{
By direct computation it is possible to check that the trigonometric polynomial $\mathbf{p}_{_{L_{2}}}$ verifies
\begin{itemize}
 \item 
$\left.\mathbf{p}_{_{L_{2}}}(\theta)^{H}\mathbf{p}_{_{L_{2}}}(\theta)+\mathbf{p}_{_{L_{2}}}(\theta+\pi)^{H}\mathbf{p}_{_{L_{2}}}(\theta+\pi)=\begin{bmatrix}
12 & 2\\
2 & 12
\end{bmatrix}>0. \quad  \right\}  \Rightarrow\textcolor{red}{\,(i)}$
 \item $\left.
\begin{tabular}{@{}c@{}}
$\mathbf{\textbf{p}}^{(2)}(0)\,{\rm e}_{2}=\begin{bmatrix}
2 & 2\\
2 & 2
\end{bmatrix}e_2=4\,{\rm e}_{2}$ \\
$\mathbf{\textbf{p}}^{(2)}(\pi)\,{\rm e}_{2}=\begin{bmatrix}
0 & 0\\
-2 & 2
\end{bmatrix}e_2=0\,{\rm e}_{2}$ \\
$\,\mathbf{\textbf{p}}^{(2)}(0)^H\,{\rm e}_{2}=\begin{bmatrix}
2 & 2\\
2 & 2
\end{bmatrix}e_2=4\,{\rm e}_{2}$ \\
\end{tabular}
\qquad\qquad \qquad \qquad \qquad\qquad \quad\right\}\Rightarrow\textcolor{red}{\,(ii)}$ 
				\item  $\lambda_{1}(\mathbf{\textbf{p}}^{(2)}(\theta+\pi))=0.$  $\qquad \quad \quad \qquad \qquad\qquad \qquad \qquad \qquad \qquad\left.\right\}\Rightarrow \textcolor{red}{(\ref{eq:final_bound_vcycle})}
$
			\end{itemize}
			In particular, the second, third and fourth conditions are the hypotheses of Lemma \ref{lemm:condition_on_s} needed   for the validation of $(ii)$. The fact that the  minimum eigenvalue function of $\mathbf{\textbf{p}}^{(2)}(\theta+\pi)$  is identically zero  implies the thesis of Lemma \ref{lem:convergence_vcycle}. Hence,  the  convergence and optimality of the V-cycle method is ensured when applied to the problem (\ref{FEM_problem}) using $P_{n,k}^{{2}}$. This is reflected on the fact that the  number of iterations needed to reach the convergence of the methods is constant when increasing the problem size. In Table \ref{tab:iterations} we show the results for the cases $a(x)=1, x^2+1, {\rm e}^{x}-2x$, using as  pre and post smoother 1 iteration of the Gauss-Seidel method. We highlight that also other smoothers are suitable for the optimality of the method. For instance, Lemma~\ref{lem:V} guarantees that the relaxed Richardson method can be used, provided that a preliminary study for the choice of the damping parameter is performed.			
\begin{table}
 \label{tab:iterations}
		\textcolor{red}{	\begin{tabular}{c|cc|cc|cc}			
		$N=d\cdot 2^t-1 $	& \multicolumn{2}{c|}{{$a(x)=1$}} & \multicolumn{2}{|c|}{{$a(x)=x^2+1$}} & \multicolumn{2}{|c}{{$a(x)={\rm e}^{x}-2x$}}\\
			\hline
			{$t$} & {TGM} & {V-Cycle} & {TGM} & {V-Cycle} & {TGM} & {V-Cycle} \\
		\hline
			4   & 6 & 6 & 7 & 7 & 6 & 7 \\
			5  & 6 & 6 & 7 & 7 & 6 & 7 \\
			6  & 6 & 6 & 7 & 7 & 6 & 7 \\
			7  & 6 & 6 & 7 & 7 & 6 & 7 \\
			8 & 6 & 6 & 7 & 7 & 6 & 7 \\
			9 & 6 & 6 & 7  & 7 & 6 & 7 \\
			10 & 6 & 6 & 7 &7& 6 & 7 \\
	\end{tabular}}
		\caption{ Two-grid and V-cycle iterations in 1D for $a(x)=1, x^2+1, {\rm e}^{x}-2x$, $tol=1 \times 10^{-6}$.}
\end{table}
}


\section{Extension to multi-dimensional case}\label{sect:multiD}

In the present subsection we give a possible extension of the convergence results in the multidimensional setting. First, we need to introduce the multi-index notation and define the objects of our analysis in more dimensions.

Let ${\bf n}:=(n_1,\ldots,n_{ {m}})$
 be a multi-index
in $\mathbb N^{ {m}}$ and set $N({r},\textbf{n}):={r}\prod_{i=1}^{ {m}} n_i$. In particular, we need to provide a generalized definition of the projector $P_{n,k}^r$ for the ${ {m}}-$level block-circulant matrix $A_N=\mathcal{A}_{\textbf{n} }(\mathbf{f})$ of dimension $N({r},\textbf{n})$ generated by a multilevel matrix-valued trigonometric polynomial $\mathbf{f}$. A complete discussion on the multi-index notation can be found in \cite{GS}.
Analogously to the scalar case, we want to construct the projectors from an arbitrary multilevel block-circulant matrix $\mathcal{A}_{\textbf{n}}(\mathbf{p})$, with $\mathbf{p}$ multivariate matrix-valued trigonometric polynomial. For the construction of the projector we can use a tensor product approach:  
 \begin{equation}
  P_{\textbf{n},\textbf{k}}=\mathcal{A}_{\bf n}(\mathbf{p}) \left(K_{\textbf{n},\textbf{k}}^{Even}\otimes I_{r}\right),
 \end{equation} 
  where $K_{\textbf{n},\textbf{k}}^{Even}$ is the $N(1,\textbf{n}) \times \frac{N(1,\textbf{n})}{2^{{ {m}}}}$ matrix defined by $K_{\textbf{n},\textbf{k}}^{Even}=K_{n_1,k_1}^{Even} \otimes K_{n_2,k_2}^{Even} \otimes \dots \otimes K_{n_{ {m}},k_{ {m}}}^{Even}$ and $\mathcal{A}_{\textbf{n}}(\mathbf{p})$ is a multilevel block-circulant matrix generated by $\mathbf{p}$.
The main goal is to combine the proof of Theorem \ref{th:tgmopt} with the multilevel techniques in \cite{NM}, in order to generalize  conditions \textit{(i)-(iii)} to the multilevel case.

In the $ { {m}}-$level setting, we are assuming that $\boldsymbol{\theta}_0\in[0,2\pi)^{ {m}}$ and $\bar{\jmath} \in\{1,\dots,r\}$ such that 
\begin{equation}\label{eqn:condition_on_f_multi}
\left\{\begin{array}{ll}
\lambda_j(\mathbf{f}(\boldsymbol{\theta}))=0 & \mbox{for } \boldsymbol{\theta}=\boldsymbol{\theta}_0 \mbox{ and } j=\bar{\jmath}, \\
\lambda_j(\mathbf{f}(\boldsymbol{\theta}))>0 & {\rm otherwise}.
\end{array}\right.
\end{equation}

The latter assumption means that the matrix $\mathbf{f}( \boldsymbol{\theta})$ has exactly one zero eigenvalue in $ \boldsymbol{\theta}_0$ and it is positive definite in $[0,2\pi)^{ {m}}\backslash\{ \boldsymbol{\theta}_0\}$.
Let us assume that, $q_{\bar{\jmath}}(\boldsymbol{\theta}_0)$ is the eigenvector of $\mathbf{f}(\boldsymbol{\theta}_0)$ associated with $\lambda_{\bar{\jmath}}(\mathbf{f}(\boldsymbol{\theta}_0))=0$. 
Moreover, define $\Omega(\boldsymbol{\theta})= \left\{\boldsymbol{\theta}+ \pi \boldsymbol{\eta}, \, \boldsymbol{\eta}\in \{0, 1\}^{ {m}} \right\}$.
Under these hypotheses, the multilevel extension of conditions \textit{(i)-(iii)}, which are sufficient to ensure the optimal convergence of the TGM in the multilevel case, is the following. Choose   $\mathbf{p}(\cdot)$ such that
\begin{itemize}
\item \begin{equation}\label{eqn:condition_on_p_multi}
	\sum_{\xi \in \Omega(\boldsymbol{\theta}) }\mathbf{p}(\xi)^{H}\mathbf{p}(\xi)>0, \quad\forall\, \boldsymbol{\theta} \in[0,2\pi)^{ {m}},
\end{equation}
 which implies that the trigonometric function
\begin{equation*}
\textbf{s}(\boldsymbol{\theta}) = \mathbf{p}(\boldsymbol{\theta})\left(\sum_{\xi \in \Omega(\boldsymbol{\theta}) }\mathbf{p}(\xi)^{H}\mathbf{p}(\xi)\right)^{-1}\mathbf{p}(\boldsymbol{\theta})^{H}
\end{equation*}
is well-defined for all $\boldsymbol{\theta} \in[0,2\pi)^{ {m}}$.
\item \begin{equation}\label{eqn:condition_on_s_multi}
	\textbf{s}(\boldsymbol{\theta} _0)q_{\bar{\jmath}}(\boldsymbol{\theta}_0) = q_{\bar{\jmath}}(\boldsymbol{\theta}_0).
\end{equation}
\item \begin{equation}\label{eqn:condition_on_s_f_multi}
\lim_{\boldsymbol{\theta} \rightarrow \boldsymbol{\theta}_0} \lambda_{\bar{\jmath}}(\mathbf{f}(\boldsymbol{\theta}))^{-1}(1-\lambda_{\bar{\jmath}}(\textbf{s}(\boldsymbol{\theta})))=c,
\end{equation}
 where $c\in\mathbb{R}$ is a constant. 
\end{itemize}
In the following we want to construct a multilevel projector  $P_{\textbf{n},\textbf{k}}$ such that the conditions (\ref{eqn:condition_on_p_multi})-(\ref{eqn:condition_on_s_f_multi}) are satisfied and, then, the optimal convergence of the TGM, applied to the problem (\ref{FEM_problem}) in the multidimensional setting, is ensured.
 In particular, starting from matrix-valued trigonometric polynomials  $\mathbf{p}_{r_\ell},$ $\ell=1,\dots, { {m}}$, we aim at defining a multivariate  polynomial $\mathbf{p}^{({ {m}})}_{{\textcolor{black}{\mathbf{r}}}}$ associated to the multilevel projector  $P_{{\textbf{n}},{\textbf{k}}}$ such that the conditions (\ref{eqn:condition_on_p_multi})-(\ref{eqn:condition_on_s_f_multi}) are satisfied.

In the following lemmas, we show that the aforementioned goal is achieved,  if  we choose the multivariate matrix-valued trigonometric polynomial
\begin{equation}\label{eqn:def_kron_p}
\mathbf{p}^{({ {m}})}_{{\textcolor{black}{\mathbf{r}}}}(\theta_1,\theta_2,\dots,\theta_{{ {m}}})=\bigotimes_{\ell=1}^{{ {m}}}\mathbf{p}_{{r_\ell}}(\theta_\ell),
\end{equation}
where $\mathbf{p}_{{\textcolor{black}{r_{\ell}}}}(\theta_\ell)\in \C^{{\textcolor{black}{r_{\ell}}}\times {\textcolor{black}{r_{\ell}}}}$ are polynomials that satisfy  conditions \textit{(i)-(iii)}.

\begin{lemma}\label{lemma:prop_kron_p}
Let $\mathbf{p}^{({ {m}})}_{{\textcolor{black}{\mathbf{r}}}}(\theta_1,\theta_2,\dots,\theta_{{ {m}}})$ be defined as in (\ref{eqn:def_kron_p}). Then,
\begin{equation*}
\sum_{\xi \in \Omega(\boldsymbol{\theta}) }\mathbf{p}^{({ {m}})}_{{\textcolor{black}{\mathbf{r}}}}(\xi)^{H}\mathbf{p}^{({ {m}})}_{{\textcolor{black}{\mathbf{r}}}}(\xi)=\bigotimes_{\ell=1}^{{ {m}}}\left(\mathbf{p}_{{r_\ell}}(\theta_\ell)^{H}\mathbf{p}_{{r_\ell}}(\theta_\ell)+\mathbf{p}_{{r_\ell}}(\theta_\ell+\pi)^{H}\mathbf{p}_{{r_\ell}}(\theta_\ell+\pi) \right).
\end{equation*}

\begin{proof}
By definition, $\mathbf{p}^{({ {m}})}_{{\textcolor{black}{\mathbf{r}}}}(\boldsymbol{\theta})=\bigotimes_{\ell=1}^{{ {m}}}\mathbf{p}_{{r_\ell}}(\theta_\ell)$, then 
\begin{small}
\begin{equation*}
\begin{split}
\sum_{\xi \in \Omega(\boldsymbol{\theta}) }\mathbf{p}^{({ {m}})}_{{\textcolor{black}{\mathbf{r}}}}(\xi)^{H}\mathbf{p}^{({ {m}})}_{{\textcolor{black}{\mathbf{r}}}}(\xi)&=\sum_{\xi \in \Omega(\boldsymbol{\theta}) } \left(\bigotimes_{\ell=1}^{{ {m}}}\mathbf{p}_{{r_\ell}}(\xi_\ell)^H\right)\left(\bigotimes_{\ell=1}^{{ {m}}}\mathbf{p}_{{r_\ell}}(\xi_\ell)\right)\\
&=\sum_{\xi \in \Omega(\boldsymbol{\theta}) } \left(\bigotimes_{\ell=1}^{{ {m}}}\left(\mathbf{p}_{{r_\ell}}(\xi_\ell)^H\mathbf{p}_{{\textcolor{black}{r_{\ell}}}}(\xi_\ell)\right)\right).
\end{split}
\end{equation*}
\end{small}
The proof is then concluded once we prove by induction on ${ {m}}$ the following equality
\begin{equation}\label{eqn:induction_p}
\sum_{\xi \in \Omega(\boldsymbol{\theta}) } \left(\bigotimes_{\ell=1}^{{ {m}}}\left(\mathbf{p}_{{r_\ell}}(\xi_\ell)^H\mathbf{p}_{{\textcolor{black}{r_{\ell}}}}(\xi_\ell)\right)\right)= \bigotimes_{\ell=1}^{{ {m}}}\left(\mathbf{p}_{{r_\ell}}(\theta_\ell)^{H}\mathbf{p}_{{r_\ell}}(\theta_\ell)+\mathbf{p}_{{r_\ell}}(\theta_\ell+\pi)^{H}\mathbf{p}_{{r_\ell}}(\theta_\ell+\pi) \right).
\end{equation}
The equation above is clearly verified for ${ {m}}=1$, indeed, by definition
\begin{gather*}
\sum_{\xi \in \Omega(\boldsymbol{\theta}) } \left(\bigotimes_{\ell=1}^{1}\left(\mathbf{p}_{{r_\ell}}(\xi_\ell)^H\mathbf{p}_{{\textcolor{black}{r_{\ell}}}}(\xi_\ell)\right)\right)= \sum_{\xi \in \{\theta_1,\theta_1+\pi \} } \left(  \mathbf{p}_{{r_1}}(\xi_1)^H\mathbf{p}_{{r_1}}(\xi_1)\right)=\\
\mathbf{p}_{r_1}(\theta_1)^{H}\mathbf{p}_{{r_1}}(\theta_1)+\mathbf{p}_{{r_1}  }(\theta_1+\pi)^{H}\mathbf{p}_{{r_1}  }(\theta_1+\pi) =\\
 \bigotimes_{\ell=1}^{1}\left(\mathbf{p}_{{r_\ell}}(\theta_\ell)^{H}\mathbf{p}_{{r_\ell}}(\theta_\ell)+\mathbf{p}_{{r_\ell}}(\theta_\ell+\pi)^{H}\mathbf{p}_{{r_\ell}}(\theta_\ell+\pi) \right).
\end{gather*}
Let us assume that equality (\ref{eqn:induction_p}) is true for ${ {m}}-1$. We have that
\begin{equation*}
\begin{split}
& \bigotimes_{\ell=1}^{{ {m}}}\left(\mathbf{p}_{{r_\ell}}(\theta_\ell)^{H}\mathbf{p}_{{r_\ell}}(\theta_\ell)+\mathbf{p}_{{r_\ell}}(\theta_\ell+\pi)^{H}\mathbf{p}_{{r_\ell}}(\theta_\ell+\pi) \right)=\\
&\left[\bigotimes_{\ell=1}^{{ {m}}-1}\left(\mathbf{p}_{{r_\ell}}(\theta_\ell)^{H}\mathbf{p}_{{r_\ell}}(\theta_\ell)+\mathbf{p}_{{r_\ell}}(\theta_\ell+\pi)^{H}\mathbf{p}_{{r_\ell}}(\theta_\ell+\pi) \right)\right]\otimes  \\
&\left(\mathbf{p}_{{r_{ {m}}}}(\theta_{ {m}})^{H}\mathbf{p}_{{r_{ {m}}}}(\theta_{ {m}})+\mathbf{p}_{{r_{ {m}}}  }(\theta_{ {m}}+\pi)^{H}\mathbf{p}_{{r_{ {m}}}  }(\theta_{ {m}}+\pi) \right)
\end{split}
\end{equation*}
The left-hand side of the latter term is a function of ${ {m}}-1$ variables $(\theta_1,\theta_2, \dots, \theta_{{ {m}}-1})$. Then, by the inductive hypothesis and from the properties of the tensor product we have
\begin{small}
\begin{gather*}
\left[\bigotimes_{\ell=1}^{{ {m}}-1}\left(\mathbf{p}_{{r_\ell}}(\theta_\ell)^{H}\mathbf{p}_{{r_\ell}}(\theta_\ell)+\mathbf{p}_{{r_\ell}}(\theta_\ell+\pi)^{H}\mathbf{p}_{{r_\ell}}(\theta_\ell+\pi) \right)\right]\otimes \\
 \left(\mathbf{p}_{{r_{ {m}}}}(\theta_{ {m}})^{H}\mathbf{p}_{{r_{ {m}}}}(\theta_{ {m}})+\mathbf{p}_{{r_{ {m}}}  }(\theta_{ {m}}+\pi)^{H}\mathbf{p}_{{r_{ {m}}}  }(\theta_{ {m}}+\pi) \right)=\\
 \left( \sum_{\substack{(\xi_1,\xi_2,\dots,\xi_{{ {m}}-1})\\ \in \\ \Omega(\theta_1,\theta_2,\dots,\theta_{{ {m}}-1})}}\bigotimes_{\ell=1}^{{ {m}}-1}\mathbf{p}_{{r_\ell}}(\xi_\ell)^{H}\mathbf{p}_{{r_\ell}}(\xi_\ell)\right)\otimes\\
 \left(\mathbf{p}_{{r_{ {m}}}}(\theta_{ {m}})^{H}\mathbf{p}_{{r_{ {m}}}}(\theta_{ {m}})+\mathbf{p}_{{r_{ {m}}}  }(\theta_{ {m}}+\pi)^{H}\mathbf{p}_{{r_{ {m}}}  }(\theta_{ {m}}+\pi) \right)=\\
 \sum_{\substack{(\xi_1,\xi_2,\dots,\xi_{{ {m}}-1})\\ \in \\ \Omega(\theta_1,\theta_2,\dots,\theta_{{ {m}}-1})}}\!\!\left[\!\left(\bigotimes_{\ell=1}^{{ {m}}-1}\mathbf{p}_{{r_\ell}}(\xi_\ell)^{H}\mathbf{p}_{{r_\ell}}(\xi_\ell)\!\right)\!\!\otimes\!\! \left(\mathbf{p}_{{r_{ {m}}}}(\theta_{ {m}})^{H}\mathbf{p}_{{r_{ {m}}}}(\theta_{ {m}})\!+\!\mathbf{p}_{{r_{ {m}}}  }(\theta_{ {m}}\!+\!\pi)^{H}\mathbf{p}_{{r_{ {m}}}  }(\theta_{ {m}}\!+\!\pi)\!\right)\!\right]\!\!=\\
 \sum_{\substack{\xi \in \{(\theta_1+l_1\pi,\dots,\theta_{{ {m}}-1}+l_{{ {m}}-1}\pi\},\\  \textbf{l}\in \{0,1\}^{{ {m}}-1}}}\left[\left(\bigotimes_{\ell=1}^{{ {m}}-1}\mathbf{p}_{{r_\ell}}(\xi_\ell)^{H}\mathbf{p}_{{r_\ell}}(\xi_\ell)\right)\otimes \mathbf{p}_{{r_{ {m}}}}(\theta_{ {m}})^{H}\mathbf{p}_{{r_{ {m}}}}(\theta_{ {m}})+\right.\\
 +\left. \left(\bigotimes_{\ell=1}^{{ {m}}-1}\mathbf{p}_{{r_\ell}}(\xi_\ell)^{H}\mathbf{p}_{{r_\ell}}(\xi_\ell)\right)\otimes \mathbf{p}_{{r_{ {m}}}}(\theta_{ {m}}+\pi)^{H}\mathbf{p}_{{r_{ {m}}}}(\theta_{ {m}}+\pi)\right]= \\
\sum_{\substack{\xi \in \{(\theta_1+l_1\pi,\dots,\theta_{{ {m}}-1}+l_{{ {m}}-1}\pi,\theta_{{ {m}}})\},\\  \textbf{l}\in \{0,1\}^{{ {m}}-1}}} \bigotimes_{\ell=1}^{{ {m}}}\mathbf{p}_{{r_\ell}}(\xi_\ell)^{H}\mathbf{p}_{{r_\ell}}(\xi_\ell)+\\
 \sum_{\substack{\xi \in \{(\theta_1+l_1\pi,\dots,\theta_{{ {m}}-1}+l_{{ {m}}-1}\pi,\theta_{{ {m}}}+\pi)\},\\  \textbf{l}\in \{0,1\}^{{ {m}}-1}}}\bigotimes_{\ell=1}^{{ {m}}}\mathbf{p}_{{r_\ell}}(\xi_\ell)^{H}\mathbf{p}_{{r_\ell}}(\xi_\ell)=\\
\sum_{\xi\in \Omega(\boldsymbol{\theta})} \bigotimes_{\ell=1}^{{ {m}}}\mathbf{p}_{{r_\ell}}(\xi_\ell)^{H}\mathbf{p}_{{r_\ell}}(\xi_\ell).
\end{gather*}
\end{small}
Then, relation (\ref{eqn:induction_p}) is verified for ${ {m}}$, and this concludes the proof.
\end{proof}
\end{lemma}
\begin{lemma}\label{lemma:positivity_kron_p}
Let $\mathbf{p}^{({ {m}})}_{\textcolor{black}{\mathbf{r}}}(\theta_1,\theta_2,\dots,\theta_{{ {m}}})$ defined as in (\ref{eqn:def_kron_p}) where $\mathbf{p}_{{r_\ell}}$, for every $\ell=1,\dots,{ {m}}$, is a polynomial  which verifies the positivity condition \textit{(i)}. Then, $\mathbf{p}^{({ {m}})}_{{r}}$ is such that the positivity condition in the multilevel setting (\ref{eqn:condition_on_p_multi}) is satisfied.
\begin{proof}
The thesis is consequence of Lemma \ref{lemma:prop_kron_p} and the matrix tensor product properties. Indeed, the eigenvalues  of a tensor product of matrices are the product of the eigenvalues of the matrices. Then, condition (\ref{eqn:condition_on_p_multi}) is trivially implied from the fact that
\begin{equation*}
\sum_{\xi \in \Omega(\boldsymbol{\theta}) }\mathbf{p}^{({ {m}})}_{{\textcolor{black}{\mathbf{r}}}}(\xi)^{H}\mathbf{p}^{({ {m}})}_{{\textcolor{black}{\mathbf{r}}}}(\xi)=\bigotimes_{\ell=1}^{{ {m}}}\left(\mathbf{p}_{{r_\ell}}(\theta_\ell)^{H}\mathbf{p}_{{r_\ell}}(\theta_\ell)+\mathbf{p}_{{r_\ell}}(\theta_\ell+\pi)^{H}\mathbf{p}_{{r_\ell}}(\theta_\ell+\pi) \right),
\end{equation*}
 and from the positivity condition \textcolor{red}{holding} in the unilevel case.
\end{proof}
\end{lemma} 
\begin{lemma}\label{lemma:tensor_s}
Let $\mathbf{p}^{({ {m}})}_{{\textcolor{black}{\mathbf{r}}}}(\theta_1,\theta_2,\dots,\theta_{{ {m}}})$ be defined as in (\ref{eqn:def_kron_p}) and it verifies (\ref{eqn:condition_on_p_multi}). Then, 
the trigonometric function
\begin{equation*}
\textbf{s}(\boldsymbol{\theta}) = \mathbf{p}^{({ {m}})}_{{\textcolor{black}{\mathbf{r}}}}(\boldsymbol{\theta})\left(\sum_{\xi \in \Omega(\boldsymbol{\theta}) }\mathbf{p}^{({ {m}})}_{{\textcolor{black}{\mathbf{r}}}}(\xi)^{H}\mathbf{p}^{({ {m}})}_{{\textcolor{black}{\mathbf{r}}}}(\xi)\right)^{-1}\mathbf{p}^{({ {m}})}_{{\textcolor{black}{\mathbf{r}}}}(\boldsymbol{\theta})^{H}
\end{equation*}
is well-defined for all $\boldsymbol{\theta} \in[0,2\pi)^{ {m}}$. Moreover, it holds 
\begin{equation}
\textbf{s}(\boldsymbol{\theta})=\bigotimes_{\ell=1}^{{ {m}}}  \textbf{s}_{{r_{\ell}}}(\theta_\ell),
\end{equation}
where  $ \textbf{s}_{{r_{\ell}}}(\theta_\ell)=\mathbf{p}_{{r_\ell}}(\theta_\ell)\left(\mathbf{p}_{{r_\ell}}(\theta_\ell)^{H}\mathbf{p}_{{r_\ell}}(\theta_\ell)+\mathbf{p}_{{r_\ell}}(\theta_\ell+\pi)^{H}\mathbf{p}_{{r_\ell}}(\theta_\ell+\pi)\right)^{-1}\mathbf{p}_{{r_\ell}}(\theta_\ell)^{H}$, for $\ell=1,\dots,{ {m}}$.
\begin{proof}
From Lemma \ref{lemma:positivity_kron_p}, we have that $\textbf{s}(\boldsymbol{\theta})$ is well-defined for all $\boldsymbol{\theta} \in[0,2\pi)^{ {m}}$. From Lemma \ref{lemma:prop_kron_p} and the properties of the tensor product, we have
\begin{equation}
\begin{split}
&\textbf{s}(\boldsymbol{\theta}) = \mathbf{p}^{({ {m}})}_{{\textcolor{black}{\mathbf{r}}}}(\boldsymbol{\theta})\left(\sum_{\xi \in \Omega(\boldsymbol{\theta}) }\mathbf{p}^{({ {m}})}_{{\textcolor{black}{\mathbf{r}}}}(\xi)^{H}\mathbf{p}^{({ {m}})}_{{\textcolor{black}{\mathbf{r}}}}(\xi)\right)^{-1}\mathbf{p}^{({ {m}})}_{{\textcolor{black}{\mathbf{r}}}}(\boldsymbol{\theta})^{H}=\\
&\bigotimes_{\ell=1}^{{ {m}}}  \mathbf{p}_{{r_{\ell}}}(\theta_\ell)\left(\bigotimes_{\ell=1}^{{ {m}}}\left[\mathbf{p}_{{r_\ell}}(\theta_\ell)^{H}\mathbf{p}_{{r_\ell}}(\theta_\ell)+\mathbf{p}_{{r_\ell}}(\theta_\ell+\pi)^{H}\mathbf{p}_{{r_\ell}}(\theta_\ell+\pi)\right]^{-1}  \right)\bigotimes_{\ell=1}^{{ {m}}}  \mathbf{p}_{{r_{\ell}}}(\theta_\ell)^H=\\
&\bigotimes_{\ell=1}^{{ {m}}}\left( \mathbf{p}_{{r_\ell}}(\theta_\ell)\left[\mathbf{p}_{{r_\ell}}(\theta_\ell)^{H}\mathbf{p}_{{r_\ell}}(\theta_\ell)+\mathbf{p}_{{r_\ell}}(\theta_\ell+\pi)^{H}\mathbf{p}_{{r_\ell}}(\theta_\ell+\pi)\right]^{-1}\mathbf{p}_{{r_\ell}}(\theta_\ell)^{H}
\right)=\bigotimes_{\ell=1}^{{ {m}}}  \textbf{s}_{{r_{\ell}}}(\theta_\ell).\end{split}
\end{equation}

\textcolor{white}{.}
\end{proof}
\end{lemma}
\begin{lemma}\label{lemma:eigenvector_one_multi}
Let $\mathbf{p}^{({ {m}})}_{_{\textcolor{black}{\mathbf{r}}}}(\theta_1,\theta_2,\dots,\theta_{{ {m}}})$ be defined as in (\ref{eqn:def_kron_p}), such that, for all $\ell=1,\dots,{ {m}}$,  
 $\mathbf{p}_{{\textcolor{black}{r_{\ell}}}}(\theta_\ell)\in \C^{{\textcolor{black}{r_{\ell}}}\times {\textcolor{black}{r_{\ell}}}}$ is a polynomial {that satisfies conditions \textit{{(i)}-{(iii)}}}. 
Let ${q}_{\mathbf{r}}=\bigotimes_{\ell=1,\dots,{ {m}}}{\rm q}_{r_\ell}$, where ${\rm q}_{r_\ell}$ is the column vector of length $r_{\ell}$ such that $ \textbf{s}_{{r_{\ell}}}( {\theta_0^{(\ell)}}){\rm q}_{r_\ell}={\rm q}_{r_\ell} $, $\ell=1,\dots,	{ {m}}$.
Then, 
\[\textbf{s}(\boldsymbol{\theta_0}){q}_{\mathbf{r}}={q}_{\mathbf{r}}, \qquad \mbox{where } \boldsymbol{\theta}_0=\left(\theta_0^{(1)},\dots,\theta_0^{({ {m}})}\right).\]
\begin{proof}
 From Lemma \ref{lemma:tensor_s}, we have that $\textbf{s}(\boldsymbol{\theta_0})=\bigotimes_{\ell=1}^{{ {m}}}  \textbf{s}_{{r_{\ell}}}(\theta_0),$ then, by definition and from the properties of the tensor product, it holds
 \begin{equation}
 \textbf{s}(\boldsymbol{\theta_0}){q}_{\mathbf{r}}=\left(\bigotimes_{\ell=1}^{{ {m}}}  \textbf{s}_{{r_{\ell}}}(\theta_0)\right)\left(\bigotimes_{\ell=1}^{{ {m}}}{\rm q}_{r_\ell}\right)=\bigotimes_{\ell=1}^{{ {m}}} \left( \textbf{s}_{{r_{\ell}}}(\theta_0){\rm q}_{r_\ell}\right)=\bigotimes_{\ell=1}^{{ {m}}} {\rm q}_{r_\ell}={q}_{\mathbf{r}}.
 \end{equation}
\end{proof}
\end{lemma}
{
\begin{lemma} \label{lemma:limit_s_multi}
Let $\mathbf{p}^{({ {m}})}_{{\textcolor{black}{\mathbf{r}}}}(\theta_1,\theta_2,\dots,\theta_{{ {m}}})$ be defined as in (\ref{eqn:def_kron_p}) such that verifies (\ref{eqn:condition_on_p_multi}). Consider  $
\textbf{s}(\boldsymbol{\theta})=\bigotimes_{\ell=1}^{{ {m}}}  \textbf{s}_{{r_{\ell}}}(\theta_\ell),
$ where $$ \textbf{s}_{{r_{\ell}}}(\theta_\ell)=\mathbf{p}_{{r_\ell}}(\theta_\ell)\left(\mathbf{p}_{{r_\ell}}(\theta_\ell)^{H}\mathbf{p}_{{r_\ell}}(\theta_\ell)+\mathbf{p}_{{r_\ell}}(\theta_\ell+\pi)^{H}\mathbf{p}_{{r_\ell}}(\theta_\ell+\pi)\right)^{-1}\mathbf{p}_{{r_\ell}}(\theta_\ell)^{H},$$ for $\ell=1,\dots,{ {m}}$, and they verify condition $(iii)$. 
Then, $\textbf{s}(\boldsymbol{\theta})$ satisfies condition (\ref{eqn:condition_on_s_f_multi}).
\begin{proof}
Without loss of generality, suppose that the order of the zero of $\lambda_{\bar{\jmath}}\left(\mathbf{f}({\theta_\ell})\right)$ in $\theta_0$ is $\varsigma \ge 2$ for $\ell=1,\dots,{ {m}}$, then the functions $1-\lambda_{\bar{\jmath}}\left( \textbf{s}_{{r_{\ell}}}(\theta_\ell) \right)$ have a zero in $\theta_0$ of order at least $\varsigma\in \mathbb{N}$ for all $\ell=1,\dots,{ {m}}$ by condition \textit{(iii)}. Hence,  the $(\varsigma-1)$-th derivative of $1-\lambda_{\bar{\jmath}}\left( \textbf{s}_{{r_{\ell}}}(\theta_\ell) \right)$ in $\theta_0$ is equal to zero. Then we have, for $\ell=1,\dots,{ {m}}$,
\begin{equation*}
\left.\lambda_{\bar{\jmath}}\left( \textbf{s}_{{r_{\ell}}}(\theta_\ell) \right)^{(\varsigma-1)}\right\vert_{\theta_0}=0.
\end{equation*}
The thesis follows by direct computation of the partial derivatives of $1-\lambda_{\bar{\jmath}}(\textbf{s}(\boldsymbol{\theta}))$ in $\boldsymbol{\theta}_0$, exploiting the fact that \begin{equation*}
\textbf{s}(\boldsymbol{\theta})=\bigotimes_{\ell=1}^{{ {m}}}  \textbf{s}_{{r_{\ell}}}(\theta_\ell) \qquad \mbox{and} \qquad
\lambda_{\bar{\jmath}}\left(\textbf{s}(\boldsymbol{\theta})\right)=\prod_{\ell=1}^{{ {m}}}  \lambda_{\bar{\jmath}}\left( \textbf{s}_{{r_{\ell}}}(\theta_\ell)\right).
\end{equation*}

\textcolor{white}{.}
\end{proof}
\end{lemma}
}
\section{Conclusions and Future Developments}
We derived the conditions which ensure the optimal convergence rate of both the TGM and the V-cycle method when applied to (multilevel) block-circulant and (multilevel) block-Toeplitz matrices. In particular, we focused on the case where the generating function $\mathbf{f}$ is a  matrix-valued trigonometric polynomial and we provide several simplifications for the validation of the theoretical conditions  in practical cases.  We also generalized the results for multilevel block-Toeplitz matrices.

As a final comment, we emphasize that the one of the main aims of the paper was to give a theoretical ground to the optimal multigrid convergence for block structured matrices, where optimal means with a convergence rate independent of the matrix size. Moreover, it provided analytical proofs of the effectiveness of standard projectors, largely used in classical applications \cite{Hack}. The numerical potency of the projectors treated in subsections \ref{sec:linear_interpolation}-\ref{sec:geometric_projector} has been exploited in many different settings (multilevel, variable coefficients case) with optimal results. \textcolor{red}{In addition to Table \ref{tab:iterations}}, see Tables 1--6 in \cite{FRTBS}  {and Table V.7 in \cite{PaolaPhD}}.

\section*{Acknowledgments}
The work of Marco Donatelli, Paola Ferrari, Isabella Furci is partially supported by Gruppo Nazionale per il Calcolo Scientifico (GNCS-INdAM).
\bibliographystyle{siamplain}
\bibliography{references}

\begin{thebibliography}{10}

\bibitem{AD}
{\sc A.~Aric\`o and M.~Donatelli}, {\em A v-cycle multigrid for multilevel
  matrix algebras: proof of optimality}, Numer. Math., 105 (2007),
  pp.~511--547.

\bibitem{ADS}
{\sc A.~Aric\`o, M.~Donatelli, and S.~Serra-Capizzano}, {\em V-cycle optimal
  convergence for certain (multilevel) structured linear systems}, SIAM J.
  Matrix Anal. Appl., 26 (2004), pp.~186--214.

\bibitem{Bhatia}
{\sc R.~Bhatia}, {\em Matrix analysis}, vol.~169 of Graduate Texts in
  Mathematics, Springer-Verlag, New York, 1997.

\bibitem{Brae}
{\sc D.~Braess}, {\em Finite elements}, Cambridge University Press, Cambridge,
  third~ed., 2007.
\newblock Theory, fast solvers, and applications in elasticity theory,
  Translated from the German by Larry L. Schumaker.

\bibitem{fem}
{\sc D.~Braess}, {\em Finite elements: Theory, fast solvers, and applications
  in solid mechanics}, Cambridge University Press, 2007.

\bibitem{LFA}
{\sc A.~Brandt}, {\em Rigorous quantitative analysis of multigrid, i. constant
  coefficients two-level cycle with l\_2-norm}, SIAM Journal on Numerical
  Analysis, 31 (1994), pp.~1695--1730.

\bibitem{tutorial}
{\sc W.~L. Briggs, V.~E. Henson, and S.~F. McCormick}, {\em A Multigrid
  Tutorial, Second Edition}, SIAM, second~ed., 2000.

\bibitem{CCS}
{\sc R.~H. Chan, Q.-S. Chang, and H.-W. Sun}, {\em Multigrid method for
  ill-conditioned symmetric {T}oeplitz systems}, SIAM J. Sci. Comput., 19
  (1998), pp.~516--529.

\bibitem{signal}
{\sc V.~Del~Prete, F.~Di~Benedetto, M.~Donatelli, and S.~Serra-Capizzano}, {\em
  Symbol approach in a signal-restoration problem involving block toeplitz
  matrices}, Journal of Computational and Applied Mathematics, 272 (2014),
  pp.~399--416.

\bibitem{tik}
{\sc M.~Donatelli}, {\em A multigrid for image deblurring with tikhonov
  regularization}, Numerical Linear Algebra with Applications, 12 (2005),
  pp.~715--729.

\bibitem{D10}
{\sc M.~Donatelli}, {\em An algebraic generalization of local fourier analysis
  for grid transfer operators in multigrid based on toeplitz matrices},
  Numerical Linear Algebra with Applications, 17 (2010), pp.~179--197.

\bibitem{DFFSS}
{\sc M.~Donatelli, P.~Ferrari, I.~Furci, D.~Sesana, and S.~Serra-Capizzano},
  {\em Multigrid methods for block-circulant and block-{T}oeplitz large linear
  systems: Algorithmic proposals and two-grid optimality analysis}, Numer.
  Linear Algebra Appl., e2356 (2020).

\bibitem{PaolaPhD}
{\sc P.~Ferrari}, {\em Toeplitz and block-toeplitz structures with variants:
  From the spectral analysis to preconditioning and multigrid methods using a
  symbol approach}.
\newblock Ph.D. Thesis, Insubria University, 2020.

\bibitem{FRTBS}
{\sc P.~Ferrari, R.~I. Rahla, C.~Tablino-Possio, S.~Belhaj, and
  S.~Serra-Capizzano}, {\em Multigrid for $\mathbb{Q}_{k}$ finite element
  matrices using a (block) {T}oeplitz symbol approach}, Mathematics, 8 (2020).

\bibitem{FS1}
{\sc G.~Fiorentino and S.~Serra}, {\em Multigrid methods for {T}oeplitz
  matrices}, Calcolo, 28 (1991), pp.~283--305 (1992).

\bibitem{FS2}
{\sc G.~Fiorentino and S.~Serra}, {\em Multigrid methods for symmetric positive
  definite block {T}oeplitz matrices with nonnegative generating functions},
  SIAM J. Sci. Comput., 17 (1996), pp.~1068--1081 (1996).

\bibitem{GS}
{\sc C.~Garoni and S.~Serra-Capizzano}, {\em Generalized locally {T}oeplitz
  sequences: theory and applications. {V}ol. {II}}, Springer, Cham, 2018.

\bibitem{qp}
{\sc C.~Garoni, S.~Serra-Capizzano, and D.~Sesana}, {\em Spectral analysis and
  spectral symbol of {$d$}-variate {$\Bbb{Q}_p$} {L}agrangian {FEM} stiffness
  matrices}, SIAM J. Matrix Anal. Appl., 36 (2015), pp.~1100--1128.

\bibitem{GV}
{\sc G.~H. Golub and C.~F. Van~Loan}, {\em Matrix computations}, vol.~3 of
  Johns Hopkins Series in the Mathematical Sciences, Johns Hopkins University
  Press, Baltimore, MD, 1983.

\bibitem{Hack}
{\sc W.~Hackbusch}, {\em Multigrid methods and applications}, vol.~4 of
  Springer Series in Computational Mathematics, Springer-Verlag, Berlin, 1985.

\bibitem{H90}
{\sc P.~Hemker}, {\em On the order of prolongations and restrictions in
  multigrid procedures}, Journal of Computational and Applied Mathematics, 32
  (1990), pp.~423--429.

\bibitem{HS2}
{\sc T.~Huckle and J.~Staudacher}, {\em Multigrid methods for block {T}oeplitz
  matrices with small size blocks}, BIT, 46 (2006), pp.~61--83.

\bibitem{Huck}
{\sc T.~K. Huckle}, {\em Compact fourier analysis for designing multigrid
  methods}, SIAM Journal on Scientific Computing, 31 (2008), pp.~644--666.

\bibitem{Kato}
{\sc T.~Kato}, {\em Perturbation theory for linear operators}, Springer-Verlag,
  Berlin-New York, second~ed., 1976.
\newblock Grundlehren der Mathematischen Wissenschaften, Band 132.

\bibitem{NN}
{\sc A.~Napov and Y.~Notay}, {\em When does two-grid optimality carry over to
  the {V}-cycle?}, Numer. Linear Algebra Appl., 17 (2010), pp.~273--290.

\bibitem{Nap11}
{\sc A.~Napov and Y.~Notay}, {\em Smoothing factor, order of prolongation and
  actual multigrid convergence}, Numerische Mathematik, 118 (2011),
  pp.~457--483.

\bibitem{RST}
{\sc R.~I. Rahla, S.~Serra-Capizzano, and C.~Tablino-Possio}, {\em Spectral
  analysis of $\mathbb{P}_k$ finite element matrices in the case of
  friedrichs--keller triangulations via generalized locally {T}oeplitz
  technology}, Numer. Linear Algebra Appl., 27 (2020), p.~e2302.

\bibitem{Re}
{\sc F.~Rellich}, {\em Perturbation theory of eigenvalue problems}, Assisted by
  J. Berkowitz. With a preface by Jacob T. Schwartz, Gordon and Breach Science
  Publishers, New York-London-Paris, 1969.

\bibitem{RStub}
{\sc J.~W. Ruge and K.~St\"{u}ben}, {\em Algebraic multigrid}, in Multigrid
  methods, vol.~3 of Frontiers Appl. Math., SIAM, Philadelphia, PA, 1987,
  pp.~73--130.

\bibitem{NM}
{\sc S.~Serra-Capizzano}, {\em Convergence analysis of two-grid methods for
  elliptic {T}oeplitz and {PDE}s matrix--sequences}, Numer. Math., 92 (2002),
  pp.~433--465.

\bibitem{negaLAA}
{\sc S.~Serra-Capizzano}, {\em Matrix algebra preconditioners for multilevel
  {T}oeplitz matrices are not superlinear}, Linear Algebra Appl., 343 (2002),
  pp.~303--319.

\bibitem{ST}
{\sc S.~Serra-Capizzano and C.~Tablino-Possio}, {\em Multigrid methods for
  multilevel circulant matrices}, SIAM J. Sci. Comput., 26 (2004), pp.~55--85.

\bibitem{Chan-Sun2}
{\sc H.-W. Sun, X.-Q. Jin, and Q.-S. Chang}, {\em Convergence of the multigrid
  method for ill-conditioned block {T}oeplitz systems}, BIT, 41 (2001),
  pp.~179--190.

\bibitem{Trot}
{\sc U.~Trottenberg, C.~W. Oosterlee, and A.~Sch\"{u}ller}, {\em Multigrid},
  Academic Press, Inc., San Diego, CA, 2001. With contributions by A. Brandt,
  P. Oswald and K. St\"{u}ben.

\end{thebibliography}
\end{document}